\documentclass{article}
\usepackage{arxiv}

\usepackage{amsmath}
    \numberwithin{equation}{section}
\usepackage{hyperref}
\usepackage[utf8]{inputenc} %
\usepackage[T1]{fontenc}    %
\usepackage{booktabs}       %
\usepackage{amsfonts}       %
\usepackage{nicefrac}       %
\usepackage{microtype}      %
\usepackage{mathabx}
\usepackage{bm}
\usepackage{stmaryrd}
\usepackage{xspace}
\usepackage{arydshln}  %
\usepackage{graphicx}
\usepackage{adjustbox}
\usepackage{algorithm}     %
\usepackage{algpseudocode}
\usepackage{cleveref}      %
\usepackage[title]{appendix}
\usepackage{graphicx,wrapfig,lipsum}
\usepackage{subcaption}
\usepackage{cite}
\usepackage{caption}

\usepackage{printlen} %

\usepackage{xcolor}
\usepackage{enumitem} %
\usepackage{fancyhdr}

\graphicspath{ {.} }

\renewcommand\vec{\mathbf}
\DeclareMathOperator*{\argmax}{argmax}

\crefname{algorithm}{algorithm}{algorithms}
\Crefname{algorithm}{Algorithm}{Algorithms}
\crefname{ineq}{inequality}{inequalities}
\Crefname{ineq}{Inequality}{Inequalities}

\newcommand{\email}[1]{\href{mailto:#1}{#1}}

\date{}
\hypersetup{
pdftitle={Data-driven Linear Solver Selection and Performance Tuning for Multiphysics Simulations in Porous Media},
pdfsubject={Data-driven Linear Solver Selection and Performance Tuning for Multiphysics Simulations in Porous Media},
pdfauthor={Yury Zabegaev, Inga Berre, Eirik Keilegavlen},
pdfkeywords={Linear solver selection,
Preconditioners,
Machine Learning,
Multiphysics,
Flow and heat transport,
Thermo-hydromechanics}
}

\begin{document}

\author{
\href{https://orcid.org/0009-0006-7095-3044}{\includegraphics[scale=0.06]{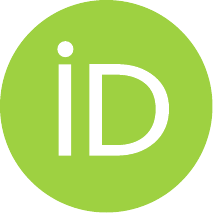}Yury Zabegaev}
\thanks{Center for Modeling of Coupled Subsurface Dynamics, Department of Mathematics, University of Bergen, Bergen, Norway}\\
\email{yury.zabegaev@uib.no}
\and
\href{https://orcid.org/0000-0002-0212-7959}{\includegraphics[scale=0.06]{orcid.pdf}Inga Berre}\footnotemark[1]\\
\email{inga.berre@uib.no}
\and
\href{https://orcid.org/0000-0002-0333-9507}{\includegraphics[scale=0.06]{orcid.pdf}Eirik Keilegavlen}\footnotemark[1]\\
\email{eirik.keilegavlen@uib.no}
}

\title{Data-driven Linear Solver Selection and Performance Tuning for Multiphysics Simulations in Porous Media}

\renewcommand{\shorttitle}{Data-driven Linear Solver Selection for Multiphysics Simulations in Porous Media}

\newpage

\maketitle

\begin{abstract}
Modeling multiphysics processes in porous media requires preconditioned iterative linear solvers to enable efficient simulations at industry-relevant scales. These solvers are typically composed of sub-algorithms that target individual physical processes. Various options are available for each algorithm, with the corresponding ranges of numerical parameters. The choices of sub-algorithms and their parameters significantly affects simulation performance and robustness. Optimizing these choices for each simulation is challenging due to the vast number of possible combinations. Moreover, optimization relies on performance data from past simulations, which becomes less representative as the simulation setup changes.
This paper addresses the problem of automated selection and tuning of preconditioned linear solvers for multiphysics simulations.
The proposed solver selection algorithm collects performance data during the run of the target simulation and continuously updates a machine learning model responsible for solver selection, resulting in an adaptively refined selection policy.
The algorithm is evaluated on two time-dependent nonlinear model problems: (i) coupled fluid flow and heat transfer in porous media and (ii) thermo-poromechanics in porous media with fractures, governed by frictional contact mechanics. These experiments demonstrate that the algorithm selects efficient and robust solvers with negligible overhead and performs comparably to a reference selection policy that has full access to the performance data of prior simulations.
Our results indicate that the proposed approach effectively addresses the challenge of solver selection and tuning, providing particular value to simulation engineers and researchers, especially when expert knowledge on linear solver tuning is not readily available.
\end{abstract}


\section{Introduction}
\label{sec:introduction}

Simulating coupled multiphysics processes in porous media is computationally challenging. Models are typically described with systems of PDEs, which are discretized and linearized, resulting in large sparse linear systems that must be solved to advance the simulation in time. Solving these linear systems is a common performance bottleneck. In industry-relevant simulations, the size of these systems varies from tens of thousands to billions of rows.
Therefore, it is only feasible to solve them with an iterative solver, which typically requires tuning of its key component -- a preconditioner. A properly selected preconditioner algorithm with tuned parameters can ensure good scalability and efficiency of the linear solver, but the optimal choice is difficult to find.
Typically, optimization for a particular problem is done empirically, combining expert knowledge with trial and error through collecting and analyzing performance data.

The model problems that we consider for the linear solver selection research are the coupled fluid flow and heat transfer in porous media and the thermo-hydromechanics (THM) problem in fractured porous media, including fracture deformation governed by contact mechanics. These are time-dependent nonlinear multiphysics problems that couple several PDEs and algebraic relations.
Preconditioners for multiphysics problems, such as the model problems considered, are usually based on sub-algorithms, that is, preconditioners originally developed for the individual single-physics components of the problem \cite{white_block-partitioned_2016,roy_block_2019}.
Extensive research has produced multiple preconditioners for the subproblems involved in these multiphysics problems and their combinations \cite{white_two-stage_2019,roy_constrained_2020,cremon_constrained_2024,bui_scalable_2020,zabegaev2025block}.
Because these sub-algorithms need to be combined, the number of possible preconditioner configurations grows exponentially compared to the optimization space of each individual single-physics problem.

Significant effort has gone into designing automatic solver selection algorithms to ensure efficient and robust numerical modeling \cite{Gries2025,Mishev2009,Huang2023,antonietti_accelerating_2023,Silva2021,Clees2010_alphaSAMG,liu2021gptune,Roy2021Bliss,zabegaev_efficient_2025}. These algorithms may be heuristic-based, such as genetic algorithms \cite{Gries2025,Mishev2009,Clees2010_alphaSAMG}, or data-driven \cite{Huang2023,antonietti_accelerating_2023,Silva2021}. The latter typically relies on collecting a large performance dataset sufficient to configure and adapt a solver for a given simulation within the considered class.
With the increasing number of configurations, the solver selection problem demands a larger performance dataset for optimization.
Collecting a single data point involves solving a linear system with a specific solver configuration and recording its performance. 
In our target applications, this will typically take from seconds to minutes, depending on the problem size and how optimal is the chosen configuration. 
Due to the combinations of the sub-algorithms configurations, the number of possible linear solver configurations can easily be in the thousands; for the simulations reported herein, the number of configurations is on the order of $10^4$. Therefore, collecting a large dataset to accurately identify the optimal solver configuration 
is justified only if 
it is used extensively in simulations over time.

In simulating highly nonlinear coupled multiphysics, it is impossible in many situations to collect an extensive dataset representative of future simulations. As the solvers behave significantly different on variants of the simulation setup, for example, based on changes in parameters, domain geometries and/or initial and boundary conditions, the collected information from a selection of simulation setups may have little relevance for other setups. In this case, we can benefit from automated solver selection during the advancement of the simulation forward in time, exploiting the similarity of the problem from one time step to another. This leads to the research question:
How can we design a selection algorithm of preconditioned linear solvers that can ensure strong simulation performance, while minimizing the time spent on exploring suboptimal solvers?

Despite a large number of possible solver configurations to select from, the available choices should not be considered arbitrary.
The possible options for sub-algorithms are typically informed by the physics of the system, so we expect the majority of the configurations to yield acceptable performance.
Although it is infeasible to always make optimal decisions over $10^4$ options and without the large performance dataset, even a small improvement over a baseline configuration is valuable in the long run, as the simulations last for days, and expert knowledge may not be readily available to select a linear solver for each particular simulation.
Therefore, we set a two-fold goal for our solver selection algorithm: First, to avoid making bad decisions, that is, picking a solver that converges slowly or fails to converge. Second, to choose an efficient solver among those that perform well.

In this paper, we propose a preconditioned linear solver selection algorithm for time-dependent multiphysics problems. This extends our earlier solver selection approach \cite{zabegaev_automated_2024} and adapts it to handle thousands of solver configurations, as opposed to a limited set.
The algorithm collects the performance data during the run of the target simulation, and utilizes the continuously updated machine learning model to guide the solver selection process. We evaluate the solver selection algorithm on the model problems designed to resemble relevant challenging simulation scenarios. We study how well the optimization process performs with incomplete performance data, aiming to balance the time spent on exploration against the benefits of larger datasets.

The manuscript is structured as follows: In \Cref{sec:model_problem}, we overview the model problems and the considered linear solver configurations for them. \Cref{sec:optimization} defines the solver selection problem, which is addressed with machine learning in \Cref{sec:ml}, where we describe the solver selection algorithm.
Simulation setups of the model problems are described in \Cref{sec:target_simulations}, and the numerical experiments with them are described in \Cref{sec:experiments}. The concluding remarks are given in \Cref{sec:conclusion}.

\section{Model Problems}
\label{sec:model_problem}
To study selection and tuning of linear solvers, we consider two model problems: First, coupled fluid flow and heat transfer in porous media, and second, THM problems in fractured porous media.
This section outlines the model problems and the considered linear solver and preconditioner configurations for them. 

\subsection{Model Problem A: Coupled Flow and Heat Transfer}

\label{sec:flow_and_transport}

Coupled flow and heat transfer can be modeled by two partial differential equations which describe the fluid mass balance and the energy balance.
A brief summary of the model is given below, for the full model we refer to \cite{coussy2004poromechanics}, while the discretization approach taken in this work is presented in \cite{porepy2024}.

The fluid flux is determined by Darcy's law:
\begin{equation}
\label{eq:darcy}
    \vec{v} = - \dfrac{\mathbf{K}}{\mu} \nabla (p - \rho \vec{g}).
\end{equation}
Here, $\mathbf{K}$ is the permeability tensor, $\mu$ is the constant fluid viscosity, $p$ is the pressure, $\vec{g}$ is the gravity vector, and $\rho$ is the fluid density.
\begin{subequations}
The fluid mass balance is given by
\begin{equation}
\label{eq:mass_balance}
    \dfrac{\partial}{\partial t} (\varphi \rho) 
    + \nabla \cdot \left( \rho \vec{v} \right) 
    = \psi,
\end{equation}
where $\varphi$ is the porosity, and $\psi$ is the fluid mass source term.
The energy balance is governed by the equation
\begin{equation}
\label{eq:energy_balance}
    \dfrac{\partial}{\partial t} (\varphi U)
    + \nabla \cdot \left( \rho c_f T \vec{v} \right) - \nabla \cdot (\kappa \nabla T)  = \psi_T,
\end{equation}
where $U$ is the porous medium internal energy, $c_f$ is the fluid specific heat capacity, $T$ is the temperature, $\kappa$ is the porous medium heat conductivity, $\psi_T$ is the energy source.
\end{subequations}

\Cref{eq:mass_balance,eq:energy_balance} are coupled through constitutive laws: The porosity and the fluid density are the functions of pressure and temperature: $\varphi = \varphi(p, T)$ and $\rho = \rho(p, T)$. Here, we do not provide the specific constitutive laws but instead refer the reader to \cite{coussy2004poromechanics}.

\subsection{Model Problem B: Thermo-poromechanics and Contact Mechanics of Fractures}
\label{sec:thm_model}
Adding a third PDE for momentum balance extends the coupled flow and heat transfer problem to the THM problem. The new independent variable introduced is the displacement, $u$.
The momentum balance equation is given by
\begin{equation}
\label{eq:momentum_balance}
    -\nabla \cdot \sigma = \psi_u,
\end{equation}
where $\sigma = \sigma(p, u, T)$ is the thermo-poroelastic stress tensor, and $\psi_u$ corresponds to the body forces. The acceleration term in the momentum balance equation is set to zero, as the quasi-static equilibrium is assumed.
The constitutive law for the porosity is extended to make it displacement-dependent: $\varphi = \varphi(p, u, T)$.

The THM problem is extended further to account for the presence of fractures, which have a significant impact on the problem physics, as they can form connected pathways that can transport substantial mass and energy through low-permeability rocks. Fractures can also deform under a combination of mechanical, hydraulic, and thermal effects. Such deformation is associated with alterations of the fracture permeability and induced seismic events \cite{ellsworth2013injection}. 
Therefore, fractures must be included in numerical models, but accurately representing them is computationally challenging due to the large difference between the porous medium grid scale and the fracture width, as well as the high aspect ratio of fracture length to width.
The problem we consider if fully described in \cite{stefansson_fully_2021}, and we provide its short overview here. Further we denote this model problem as the \textit{contact-THM} problem.

The fractures are inserted into the model THM problem with the Discrete Fracture Modeling (DFM) approach, which treats the fractures as lower-dimensional geometric objects, resulting in a mixed-dimensional problem formulation, as introduced in \cite{alboin1999domain,alboin2002modeling}. 
The mass and energy balance is imposed on fractures by equations similar to \eqref{eq:mass_balance}  and \eqref{eq:energy_balance}.
The force balance equation is imposed on the fracture sides in the form of Newton's 3rd law. Additionally, the contact mechanics inequalities are imposed: the non-penetration inequalities and the Coulomb friction law, as defined in \cite{kikuchi1988contact}. Together, they describe three distinct states the fracture can be in: sticking, sliding or opening of the fracture sides, which alters the resulting matrix properties for the corresponding fractures. Note that the contact problem does not directly affect solver selection, as we consider only one sub-algorithm to address it. Therefore, we do not discuss it in detail here. For further information, the reader is referred to \cite{stefansson_fully_2021}.

\subsection{Linear Solvers for the Model Problems}
\label{sec:linear_solvers_general}

The considered model problems are nonlinear, so the solution strategy is based on the Newton method. In each of its iterations, the Jacobian-based linear system must be solved with a preconditioned iterative linear solver. In this paper, we consider the generalized minimal residual (GMRES) method, as the Jacobian matrix is non-symmetric. 

We describe the preconditioners for the contact-THM problem below, and the preconditioners for the coupled mass and heat transfer problem can be derived by dropping the treatment of the momentum balance and contact mechanics.
The Jacobian matrix of the contact-THM problem consists of three main blocks, corresponding to the momentum, mass and energy balance equations \eqref{eq:momentum_balance}, \eqref{eq:mass_balance}, and \eqref{eq:energy_balance}:
\begin{equation}
\label{eq:jac}
\begin{array}{ccc}
J = &
\left[\ \begin{matrix}
    J_{uu} & J_{up} & J_{uT} \\
    J_{pu} & J_{pp} & J_{pT} \\
    J_{Tu} & J_{Tp} & J_{TT} \\
\end{matrix}\ \right]
&
\begin{array}{l}
     \text{ mechanics}\\
     \text{ fluid flow}\\
     \text{ heat transfer}
\end{array}
\\
&
\begin{array}{ccc}
u \ \ & \ \ p \ \ & \ \ T 
\end{array}
&
\end{array}
\end{equation}
Here, we use the row and the column captions to denote the equation and the independent variable that correspond to the submatrix. For instance, the submatrix $J_{uu}$ is produced by differentiating the momentum balance equation by the displacement variable.
The flow and heat transfer model problem can be represented by a $2\times2$ block matrix by considering only the $p$ and $T$ submatrices. 
For the more detailed presentation of the block matrix corresponding to the contact-THM problem, the reader is referred to \cite{stefansson_fully_2021}.

Preconditioners for multiphysics problems, such as the THM problem, are typically based on a block factorization of the full Jacobian matrix, see \cite{benzi_numerical_2005} for details. This technique decouples the submatrices on the main diagonal so that their inverses can be approximated individually by sub-algorithms tailored for particular single-physics linear systems. We will follow this paradigm for solving \eqref{eq:jac}, and hence our design of linear solvers can be described in terms of the decoupling technique chosen and the preconditioner applied in the sub-algorithms.

The momentum balance is typically decoupled first from the flow and transport subproblem, see \cite{bui_scalable_2020,white_two-stage_2019}. The decoupling is based on the fixed-stress approximation, as described in \cite{white_block-partitioned_2016,kim_stability_2011,kim_unconditionally_2018}. The decoupled momentum-balance submatrix is approximately solved with the algebraic multigrid method (AMG) \cite{ruge_amg_1986,griebel_amg_2003}, which is suited for elliptic-type problems.

The flow and heat transfer subproblem is tightly coupled, and depending on the flow regime, the heat transfer can be dominated by the advective or diffusive term in \eqref{eq:energy_balance}. If the advective term prevails, the problem is typically addressed by the constrained pressure residual method (CPR) \cite{wallis_constrained_1985,roy_block_2019,white_two-stage_2019}.
The method consists of two stages: (i) the pressure subproblem is approximately solved with AMG, and (ii) the coupled pressure-temperature submatrix is approximately solved with the incomplete LU factorization (ILU). If the diffusive term dominates, it is efficient to apply the AMG method to the coupled pressure-temperature matrix, see, for instance \cite{roy_block_2019}. We refer to this approach as System-AMG in this paper. The challenge is that the ratio between the two terms can vary within one simulation, and it is hard to determine in advance, which method will perform better.

Extending the preconditioner for the contact-THM problem requires several additional steps to treat the contact mechanics submatrix and the submatrices of the lower-dimensional processes. Significant to this work, none of the extra steps introduce major decisions to be made to realize a particular solver, and we will therefore not discuss this further, but refer the reader to \cite{zabegaev_efficient_2025,zabegaev2025block} for details.

The main decisions in solver design are thus, first, whether to treat the mass-energy subsystem by CPR or solve them jointly, and second, how to apply AMG in the various sub-algorithms.
The latter is a non-trivial matter, as AMG should be seen not as a single algorithm, but as a class of algorithms \cite{Gries2025}. This is illustrated by two examples: First, the classical formulation of AMG \cite{rugestuben} and the smoothed aggregation AMG \cite{Vanek1996} construct the operators used in the method differently, resulting in performance differences for a given problem. 
Second, at its core, any AMG method applies inexpensive linear subsolvers, such as Jacobi or Gauss--Seidel, which are referred to as ``smoothers'' in the multigrid literature. The choice of smoothers can significantly impact performance. We will discuss AMG parameters in some more detail in Section \ref{sec:target_simulations}.

To achieve good performance of AMG, its parameters must be tuned to balance per-iteration performance against approximation accuracy, which in turn determines the number of iterations. The optimal values are problem-dependent: they differ for the mechanics, fluid, and energy AMG subsolvers and can vary between problem setups.
Considerable effort goes into tuning AMG parameters for specific problems, including the use of machine learning methods \cite{Gries2025,Huang2023,antonietti_accelerating_2023}.
In our case, the nonlinear multiphysics couplings can alter the properties of the resulting matrices and affect optimal parameter configurations. Choosing the appropriate AMG algorithms and tuning their parameters is therefore a major challenge in solver selection for multiphysics problems, though not the only one.

\section{Solver Selection Problem}
\label{sec:optimization}

In this section, we describe the solver selection problem as an optimization problem. We consider three key components of the problem: (i) linear solver configurations, (ii) the simulation state, which determines the difficulty of the current linear system, and (iii) the optimization target function.

\subsection{Solver Configuration Space}

\label{sec:solver_space}
\begin{figure}[htb]
    \centering
    \includegraphics[width=1\linewidth]{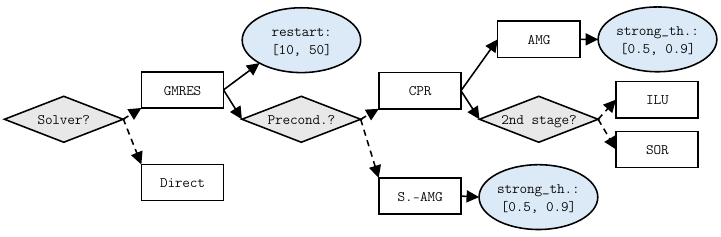}
    \caption{The flowchart represents the solver configuration space for the coupled flow and transport problem. Diamond blocks correspond to algorithmic choices, oval blocks represent numerical parameters to tune, and rectangular blocks denote algorithms to choose from. Dotted arrows indicate that a single path must be chosen, while solid arrows connect components of a single linear solver configuration.
    }
    \label{fig:solver_space}
\end{figure}

A preconditioned linear solver configuration for a multiphysics problem typically consists of multiple algorithmic decisions, as discussed in \Cref{sec:linear_solvers_general}. The range of candidate solver configurations can be visualized as a flowchart. As an example, \Cref{fig:solver_space} illustrates such a flowchart for the coupled flow and heat transfer model problem, described in \Cref{sec:flow_and_transport}. The example concerns the choice between using the preconditioned GMRES method or a direct solver, the latter being optimal for small problem sizes. The GMRES \texttt{restart} parameter can be selected within given bounds, and two preconditioners are considered: CPR and System-AMG, as described in \Cref{sec:linear_solvers_general}. The AMG \texttt{strong threshold} parameter can be tuned, while for the second CPR stage, two options are available: ILU or the cheaper SOR method. We reiterate that this example is used to illustrate the solver space, in practice, we consider significantly more choices in the numerical experiments.

\begin{figure}[ht]
    \centering
    \includegraphics[width=1\linewidth]{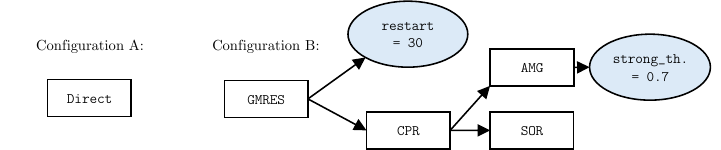}
    \caption{Two solver configurations taken from the solver configuration space illustrated in \Cref{fig:solver_space}. Configuration A represents a direct solver, while Configuration B represents GMRES, preconditioned by the CPR method.
    }
    \label{fig:chosen_solvers}
\end{figure}

The term \textit{solver configuration} is used to refer to a complete algorithm for solving a linear system, with all algorithmic decisions made and all numerical parameters specified. It is denoted by $a \in \mathcal{A}$, which is the \textit{solver configuration space}.
Two example solver configurations are illustrated in \Cref{fig:chosen_solvers}. 
Most of the linear solver tuning lies in the preconditioner, so the term ``preconditioner configuration'' could be used interchangeably. However, we prefer ``solver configuration'', as it can also include parameters outside the preconditioner, such as the GMRES restart parameter, while the tuning methodology remains the same.

We observe that $\mathcal{A}$ is not a Cartesian product of all individual algorithmic decisions, since some combinations are not practically meaningful. For instance, the choice of a specific CPR method configuration is irrelevant if the direct solver is selected in \Cref{fig:solver_space}.

Solver configurations $a \in \mathcal{A}$ need to be encoded into uniformly shaped vectors to be used as input in a machine learning model. The categorical decisions within a single configuration $a$ are encoded in a binary sequence that preserves the structural relations between decisions: One or zero denotes the presence or absence of a particular decision in the solver configuration.

The numerical decisions are similarly encoded into a vector, so that each tunable parameter has a prescribed position in it. The vectors of the encoded categorical and numerical decisions are then concatenated. The values of the numerical decisions that are not present in the encoded algorithm are filled with zeros. 
The examples of the encoded categorical and numerical choices are shown in \Cref{tab:cat_num_encoding}.

\begin{table}[ht]
\setlength{\abovecaptionskip}{5pt}  %
\centering
\begin{tabular}{lrrrrrr:rrr}
\toprule
Solver configuration & \multicolumn{6}{c:}{Categorical encoding} & \multicolumn{3}{c}{Num. encoding} \\
\midrule
Decision name &
\adjustbox{valign=b}{\rotatebox{-90}{\texttt{Direct}}} & 
\adjustbox{valign=b}{\rotatebox{-90}{\texttt{GMRES}}} & 
\adjustbox{valign=b}{\rotatebox{-90}{\texttt{S.-AMG}}} & 
\adjustbox{valign=b}{\rotatebox{-90}{\texttt{CPR}}} & 
\adjustbox{valign=b}{\rotatebox{-90}{\texttt{SOR}}} & 
\adjustbox{valign=b}{\rotatebox{-90}{\texttt{ILU}}} &
\adjustbox{valign=b}{\rotatebox{-90}{\texttt{restart}}} &
\adjustbox{valign=b}{\rotatebox{-90}{\texttt{strong\_th.\!1}}} & 
\adjustbox{valign=b}{\rotatebox{-90}{\texttt{strong\_th.\!2}}} \\
\midrule
Configuration A & \texttt{[1} & \texttt{0} & \texttt{0} & \texttt{0} & \texttt{0} & \texttt{0} & \texttt{0} & \texttt{0} & \texttt{0]} \\
Configuration B & \texttt{[0} & \texttt{1} & \texttt{0} & \texttt{1} & \texttt{1} & \texttt{0} & \texttt{30} & \texttt{0} & \texttt{0.7]} \\
Configuration C & \texttt{[0} & \texttt{1} & \texttt{1} & \texttt{0} & \texttt{0} & \texttt{0} & \texttt{50} & \texttt{0.5} & \texttt{0]} \\
\bottomrule
\end{tabular}
\caption{
Structure-preserving encoding of the categorical and numerical choices for the solver configuration space shown in \Cref{fig:solver_space}.
Two demonstrated vectors (rows of the table) correspond to the solver configurations A and B from \Cref{fig:chosen_solvers}. 
The last row corresponds to the System-AMG preconditioner, which is not shown in the figure.
Two numerical parameters, \texttt{strong\_th.\!1} and \texttt{strong\_th.\!2}, correspond to the System-AMG and the AMG in the first stage of CPR, respectively. They are treated separately, even though they represent the same parameter in two different instances of the same algorithm.
}
\label{tab:cat_num_encoding}
\end{table}

\subsection{Simulation State and Context}
\label{sec:context}

We denote the \textit{full simulation state} by $\chi$. It includes all the information about the simulation and its current state: the geometry, boundary conditions, injection rates, the current solution, etc.
Any state $\chi$ is associated with a single linear system. $\chi$ includes both simulation properties that remain constant and those that change during a simulation, such as the current temperature field in the simulation.
An ideal solver selection algorithm should account for the full simulation state, as it determines the difficulty of the resulting linear system. Moreover, different solver configurations may be optimal for different simulations, or even for different time steps within the same simulation.

In practice, $\chi$ is not available to the solver selection algorithm, as incorporating it would require an impractically large and expensive machine learning model.
Instead, $\chi$ must be compressed into a compact representation of the problem that includes only a fixed number of features. There is no single ``correct'' feature list; instead, one must be constructed for each problem using physical and mathematical insight into its specific nature.
We denote the \textit{simulation context} by $c$, a vector of $n$ real numbers that encode the selected feature list. The simulation context space is given by $\mathcal{C} \subset \mathbb{R}^n$.

Typically, the simulation context can include the Jacobian matrix size, mesh quality metrics, the CFL number, and dimensionless parameters relevant to the specific problem. For problems driven by source terms and boundary conditions, it is useful to include the magnitude of the corresponding terms.

Since the dimension of the context $c$ is much smaller than that of the full state $\chi$, some information is inevitably lost. As a result, multiple significantly different linear systems may correspond to identical contexts. Although this is unavoidable, a well-chosen context space for a given problem can minimize the negative impact of such cases on solver selection. 

\subsection{Optimization Target Function}
\label{sec:reward}
The performance tuning algorithm aims to select the optimal solver configuration $a^{\star} \in \mathcal{A}$ that solves the given linear system successfully and minimizes the time spent on solving it. The optimization problem is expressed in terms of a success indicator function $S$ and a reward function $R$. Consider a linear system and the corresponding simulation state $\chi$. The linear solver $a \in \mathcal{A}$ is selected to solve it. Then, the success and the reward functions for it are given by:
\begin{equation}
\label{eq:success_reward}
    S(a, \chi) = \begin{cases}
        \text{true } \quad \text{if solved successfully,} \\
        \text{false} \quad \text{otherwise.} \\
    \end{cases}
    \quad
    R(a,\chi) = -\text{log }T(a,\chi).
\end{equation}
Here, $T(a,\chi)$ denotes the time to construct the preconditioner and solve the linear system with $a$. The successful solution corresponds to the relative residual norm decrease up to a fixed tolerance. If the solver reaches a given iteration limit, this is treated as failure. We assume that the solution tolerance and the iteration limit are the same for all linear systems.

Formally, $R(a,\chi)$ is not a deterministic function because of noise in the time measurements. We assume this noise is normally-distributed and stationary. In this paper, we omit the randomness in the notation and treat $R(a,\chi)$ as the expectation with respect to this noise.

The optimal solver for the given linear system with the corresponding simulation state $\chi$ is defined by:
\begin{equation}
\label{eq:optimization_reward}
    a^\star = \argmax_{a \in \mathcal{A}} 
 R(a,\chi) \quad \text{s.t.} \quad S(a,\chi) = \text{true},
\end{equation}
The optimization problem is solved after each encountered linear system: within Newton iterations for a nonlinear problem and in each time step for a linear problem.

\section{Solver Selection Algorithm}
\label{sec:ml}

Machine learning is applied to the optimization problem \eqref{eq:optimization_reward}. As stated in \Cref{sec:introduction},
we select from $\sim 10^4$ available solver configurations, and
we aim to avoid spending significant time collecting the performance dataset, i.e. solving numerous linear systems with suboptimal linear solvers. Additionally, the machine learning method must not add significant overhead relative to linear solve time.

These constraints inform the design choices we make in developing the solver selection algorithm.
The proposed approach builds on our previous work \cite{zabegaev_automated_2024}, with key extensions and differences outlined at the end of this section. 
We approximate the success and the reward functions \eqref{eq:success_reward} with two machine learning models that form a sequential pipeline. We refer to them as two stages, and both stages observe the same simulation context $c$, which represents the full simulation state $\chi$: 
\begin{enumerate}[noitemsep]
    \item A \textbf{classifier} predicts which configurations are likely to succeed: $s(a,c) \approx S(a,\chi)$. The others are discarded and do not proceed to the second stage;
    \item A \textbf{regression model} estimates the expected reward of the remaining configurations: $r(a,c) \approx R(a,\chi)$.
\end{enumerate}

\begin{figure}[htb]
    \centering
    \includegraphics[width=0.8\linewidth]{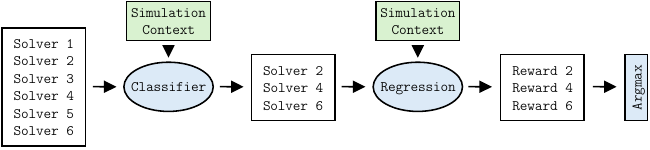}
    \caption{
    Machine learning pipeline application. In this example, the solver configuration space consists of 6 configurations. First, the classifier model discards solver configurations predicted to fail and retains only those expected to succeed. Then, the regression model predicts the rewards for those configurations. The configuration with the largest predicted reward is selected. Both the classifier and the regression models use the same solver encoding and simulation context as input.  
    }
    \label{fig:ml_pipeline}
\end{figure}

The application of a pipeline is visualized in \Cref{fig:ml_pipeline}. We specify the particular machine models used for $s$ and $r$ later in this section.
The solver configuration is selected based on the pipeline's predictions:
\begin{equation}
\label{eq:solver_selection_pipeline}
    a_\text{pred} = \argmax_{a \in \mathcal{A}} r(a,c) \quad \text{s.t.} \quad s(a,c) = \text{true}.
\end{equation}
In practice, the pipeline makes predictions for all the configurations $a \in \mathcal{A}$, and the one with the maximum predicted reward is selected.
Real-valued numerical decisions in $\mathcal{A}$ should be discretized with a specified step size.

The classifier stage is fitted with all the available performance data, including both the succeeded and failed linear solve attempts. The regression stage is fitted only with the performance data that corresponds to the successful linear solve attempts. 

The input for the machine learning pipeline consists of the solver configuration encoding $a$ and the simulation context $c$, which are described in \Cref{sec:solver_space,sec:context}, respectively.
As discussed in \Cref{sec:context}, $c$ does not fully approximate $\chi$. If the context is poorly engineered, multiple identical contexts can arise from significantly different states, leading to different predictions for the same algorithm configuration. 
Therefore, for a given context $c$, the machine learning models approximate the average reward and success expectations for all states $\chi$ that produce this context.
In the worst case, the context becomes irrelevant, so ensuring it contains descriptive features is necessary for its reliability. 

\subsection{Initial Exploration}
\label{sec:initial_exploration}

In cases where little performance data for the target simulation is available in advance, performance data collection begins at the target simulation run time. That is, we do not execute additional simulations solely to analyze the solver configuration performance. 
We refer to this process as \textit{initial exploration}. During this process, a new solver configuration from $\mathcal{A}$ is selected randomly each time, e.g. for each new Newton method iteration or time step.
If the linear solver converges, it is still switched to a new random configuration.
If it fails, the next random configuration is applied to the same linear system, and this process continues until the system is successfully solved. After the successful or failed solution attempt, the performance data is added to the dataset that includes:
\begin{itemize}[noitemsep]
    \item the chosen solver configuration encoding $a$;
    \item the simulation context $c$ based on the current simulation state $\chi$;
    \item the success or failure flag $S$;
    \item the reward $R$ produced by measuring the run time.
\end{itemize}

The number of linear systems solved in the initial exploration phase \texttt{num\_initial} is fixed.
We emphasize that setting \texttt{num\_initial} to the number of configurations in $\mathcal{A}$ does not imply that an extensive performance dataset has been collected, since the data only reflects how linear solvers perform for certain simulation states. Hence, this dataset may still be insufficient to predict the performance of all configurations across all possible states.

Minimizing the initial exploration time and applying the solver selection algorithm as early as possible can potentially save a significant share of the simulation time. However, this strategy can also be risky, as some simulation states that occur later within a simulation can be misrepresented. Additionally, some solver configurations may never be tested, potentially causing good configurations to be overlooked.

The initial exploration phase may be skipped entirely if the past simulation performance data is available and considered relevant. However, the more changes are made in the simulation setup, the less representative the historic data becomes.
Constructing a metric that would define the distance between the models, and therefore the relevance of the performance data, seems infeasible. Similarly, it is challenging to include all the possible differences between the simulations in the limited context vector for the machine learning model, described in \Cref{sec:context}. Therefore, instead of relying on the massive historical performance dataset and a heavy machine learning model to generalize across different models, we collect performance data on the fly. We believe it is a better choice given our design constraints discussed in \Cref{sec:introduction}, particularly because the current simulation's performance data is more relevant than any previously obtained data.

\subsection{Incremental Learning}
\label{sec:incremental_learning}

The issue of minimizing the initial exploration time is addressed by the \textit{incremental learning} mechanism: to continuously update the machine learning model during the simulation with the feedback to its decisions. This allows beginning the solver selection process when the machine learning approximations are relatively coarse and refine them later by the incremental updates.

The machine learning pipeline of the solver selection algorithm is integrated into the simulation time step or Newton loop and consists of two stages: (i) inference -- selecting the optimal configuration, and (ii) training -- updating the model with feedback from previous decisions. Since updating the machine learning model is expensive, feedback data is accumulated in fixed-size batches, and the model is updated only when a batch is full.

The machine learning pipeline is trained for the first time after the initial exploration phase is over. The integration of the solver selection algorithm to the simulation main loop is described as follows. At each Newton iteration or time step with the corresponding simulation state $\chi$: 
\begin{enumerate}[noitemsep]
    \item The machine learning pipeline makes predictions for all the solver configurations in $\mathcal{A}$. The one with the maximum expected reward is chosen;
    \item The chosen linear solver attempts to solve the linear system;
    \item The solution attempt data is saved into the performance dataset;
    \item If the feedback data batch is full, the machine learning pipeline is updated with the new data;
    \item If the system is solved successfully, the simulation proceeds to the next linear system; otherwise, the algorithm repeats.
\end{enumerate}

A situation when all the configurations are predicted to fail is unlikely but not impossible. This situation indicates that either: (i) there is truly no solver for this problem, (ii) the failure criterion, e.g. iterations limit is too strict, or (iii) the classifier stage of the pipeline is confused. Since these situations are indistinguishable for us, we select a random configuration. In case it solves the given problem, it may drive the classifier out of confusion. Otherwise, the simulation fails to solve the linear system, and the standard failure procedure is executed, e.g. decrease the simulation time step and retry.

In a truly incremental update, only the new data is used for the update, and the full dataset is not stored in the program memory.
The alternative is to store the full dataset and utilize it to refit the model every time the new data is received. Although this approach is significantly more expensive, it may be reasonable within our constraints. In our target application, 
a simulation requires solving on the order of $10^4$ linear systems, resulting in a corresponding number of data points. The update time of the machine learning pipeline can still be on the order of seconds, while the linear solve time may be on the order of minutes or more, depending on the problem scale.

These estimations motivate us to fully retrain the machine learning pipeline each time, since the dataset is small enough to store in memory during simulation, and the retraining time is negligible compared to the overall simulation workload. Additionally, this greatly reduces concerns about the convergence of the machine learning models. The machine learning models we use in the pipeline are the gradient boosting classification and regression models.

The key difference from our previous work \cite{zabegaev_automated_2024} is that we now operate in a setting with a much larger number of available solver configurations. We therefore adapted the method by adding a classifier machine learning model and fully adopting a greedy strategy, so the algorithm always selects the configuration with the highest predicted reward.

\section{Simulation Setups}
\label{sec:target_simulations}

The model problems described in \Cref{sec:model_problem} are discretized using Cartesian (for the problem of coupled flow and heat transfer) and simplex (for contact-THM) grids. The grid construction is handled by Gmsh \cite{gmsh}. For the mixed-dimensional grid, this approach ensures that the lower-dimensional subdomains (fractures) geometrically coincide with a set of faces in the surrounding higher-dimensional grid.  
The spatial discretization is based on a family of cell-centered finite volume schemes. The multi-point flux approximation (MPFA) method \cite{Aavatsmark2002} is applied to discretize the diffusive scalar fluxes (mass and energy), whereas the multi-point stress approximation (MPSA) method \cite{Nordbotten2016} is used to discretize the thermo-poromechanical stress.
In fractured domains, coupling of discretizations on different subdomains follows the approach defined in \cite{nordbotten2019unified}.
The first-order upwinding scheme is utilized to discretize the advective fluxes of energy and inter-dimensional fluid flux; the values of the fluid flux $\vec{v}_i$ and the inter-dimensional fluid flux $v_j$ are computed from the solution at the previous Newton iteration. The temporal discretization is performed using the implicit Euler scheme.

The discretized governing equations result in a nonlinear system that must be solved to advance the simulation in time. Due to the contact mechanics relations, the system is only semi-smooth.
To solve it, we employ a semi-smooth Newton method \cite{ito2003semi}, which involves inverting the generalized Jacobian matrix. This matrix coincides with the original Jacobian in smooth regions and adopts derivative values from one side or the other in non-smooth regions \cite{ito2003semi,berge_finite_2020}.
Implementation of the discretization and the other simulation procedures is handled by PorePy \cite{porepy2021,porepy2024}.

We present two sets of experiments, Sequence A and Sequence B, which consider, respectively, coupled fluid flow and heat transfer in a porous medium and contact thermo-hydromechanics in a fractured porous medium. Details of the sequences, including simulation parameters and the spaces of solvers made available to our selection algorithm, are given below.

In practical simulations, it is common to consider multiple similar simulations, e.g., to test sensitivity to input parameters or in optimization workflows.
Such repeated simulations can be expected to produce highly relevant data for the solver selection algorithm.
To mimic this situation, we consider a simulation workflow where multiple similar but not identical simulations are performed: For Sequence A, we vary the permeability field between individual simulations, while for Sequence B, the location of injection and production wells are changed. Linear solver performance data is shared within each sequence, but not between Sequences A and B.

\subsection*{Common Linear Solvers Description}

\begin{table}[htb]
\centering
\setlength{\abovecaptionskip}{5pt}  %
\begin{tabular}{lp{4.6cm}p{4.7cm}}
\toprule
Parameter keyword & Range of values & Clarification \\
\midrule
\texttt{hypre\_boomeramg\_strong\_threshold} & \texttt{\{0.5, 0.6, 0.7, 0.8, 0.9\}} & Threshold for deciding which matrix connections are considered ``strong'' in coarsening.\\
\texttt{hypre\_boomeramg\_agg\_nl}           & \texttt{\{0, 1, 2\}} & Number of levels of aggressive coarsening. \\
\texttt{hypre\_boomeramg\_relax\_type\_all}         & \texttt{\{symmetric-SOR/Jacobi, l1scaled-Jacobi, SOR/Jacobi, Jacobi\}} & Relaxation (smoothing) method used on all grid levels. \\
\texttt{hypre\_boomeramg\_cycle\_type}       & \texttt{\{V, W\}} & Multigrid cycling strategy. \\
\texttt{hypre\_boomeramg\_grid\_sweeps\_all} & \texttt{\{1, 2, 3\}} & Number of relaxation sweeps to perform on every grid level. \\
\hdashline
\texttt{gamg\_threshold} & \texttt{\{0, 0.01, 0.05, 0.1\}} & Threshold for deciding which matrix connections are considered ``strong'' in coarsening. \\
\texttt{gamg\_agg\_nsmooths} & \texttt{\{0, 1\}} & Number of smoothing iterations applied when forming aggregates in coarsening. \\
\texttt{gamg\_aggressive\_coarsening} & \texttt{\{1, 2\}} & Number of levels of aggressive coarsening. \\
\texttt{mg\_levels\_pc\_max\_it} & \texttt{\{sor, pbjacobi\}} & Relaxation (smoothing) method used on all grid levels. \\
\texttt{mg\_cycle\_type} & \texttt{\{V, W\}} & Multigrid cycling strategy. \\
\texttt{mg\_levels\_ksp\_max\_it} & \texttt{\{1, 2, 4\}} & Number of relaxation sweeps to perform on every grid level. \\
\bottomrule
\end{tabular}
\caption{AMG parameters considered for tuning, both for \texttt{BoomerAMG} and \texttt{GAMG} implementations. The parameter keywords are used as they are denoted in PETSc. Lower values for \texttt{hypre\_boomeramg\_strong\_threshold} and \texttt{gamg\_threshold} lead to denser coarse-grid operators, causing a potentially slower but more accurate preconditioner.
}
\label{tab:amg_params}
\end{table}

Further in this section, we discuss the preconditioned linear solver configurations considered for each sequence. Both are based on the same sub-algorithms, so we list their shared features. The preconditioners are based on PETSc \cite{petsc_user_manual} and its Python binding PETSc4Py \cite{petsc4py}. Specifically, the stationary preconditioners \texttt{Jacobi} and \texttt{SOR}, and the \texttt{ILU} algorithm are implemented in PETSc, along with their block versions, denoted such as \texttt{Block-Jacobi}. If the block version is used, the matrix is always considered in the cell-based arrangement, so that the degrees of freedom defined in the same grid cell are located contiguously in memory. The block size corresponds to the number of degrees of freedom per cell. The number \texttt{k} in \texttt{ILU(k)} corresponds to the number of levels of fill.

Sequences A and B utilize AMG preconditioner variants. Two AMG implementations are considered: \texttt{BoomerAMG} by HYPRE \cite{hypre} and \texttt{GAMG} by PETSc \cite{petsc_user_manual}. HYPRE is accessed through the PETSc interface.
\texttt{BoomerAMG} implements the classical AMG algorithm, while \texttt{GAMG} corresponds to the smoothed-aggregation AMG. Therefore, some algorithm parameters are different. We list the AMG parameters available for tuning in \Cref{tab:amg_params}.

\subsection{Sequence A: Coupled flow and heat transfer in porous media with varying permeability}
\label{sec:setup A}
The first sequence corresponds to the coupled flow and heat transfer problem, described in \Cref{sec:flow_and_transport}. We consider a 3D cuboid domain of porous media filled with water. The permeability and porosity fields are taken from the SPE10 benchmark study,
which represents a model with a challenging anisotropic and heterogeneous structure of a subsurface reservoir \cite{SPE10}. We use only the x-component of the permeability field, making it isotropic, but still heterogeneous. The permeability values are distributed within $(10^{-17}, 10^{-11})$ $\text{m}^{2}$.

The full SPE10 model consists of \texttt{(85, 220, 60)} Cartesian cells.
The simulations within Sequence A are based on different 3D slices of the permeability and porosity fields from this model:
Each slice is a cuboid of shape \texttt{(42, 100, 20)} Cartesian cells or $(121.92 \times 304.8 \times 128.016)$ $\text{m}^3$, where the third dimension is associated with the sharpest alterations in porosity and permeability. 
In total, the sequence includes 15 simulations. 

\begin{figure}[t]
    \centering
    \begin{minipage}{0.49\textwidth}
        \centering
        \includegraphics[width=\linewidth]{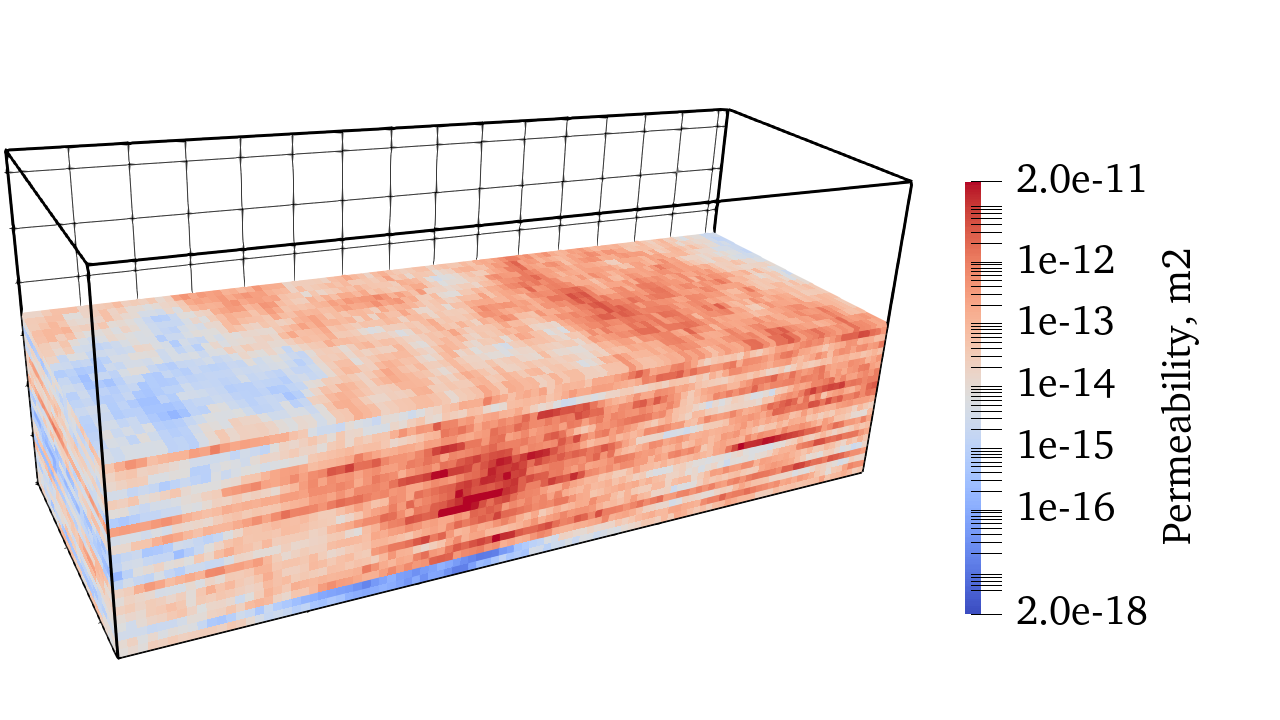}
    \end{minipage}\hfill
    \begin{minipage}{0.49\textwidth}
        \centering
        \includegraphics[width=\linewidth]{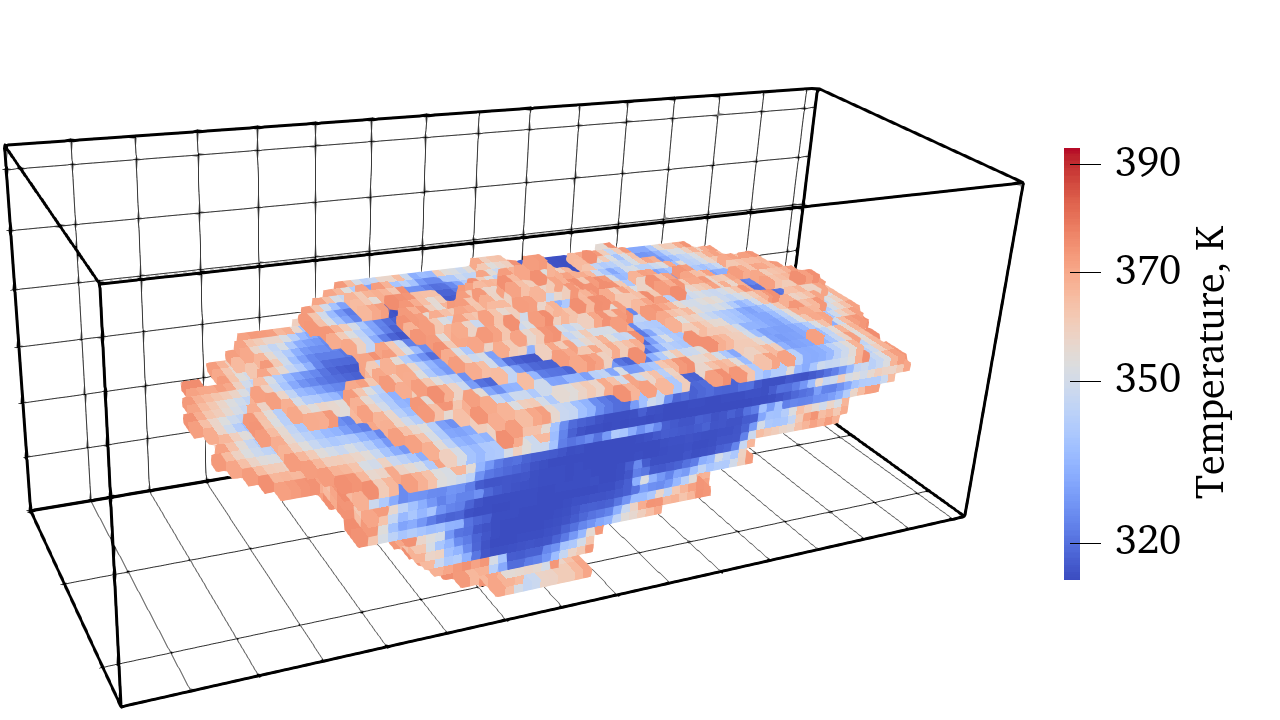}
    \end{minipage}
    \caption{One of the simulations in Sequence A, the simulation time is 150 days. Left: permeability field slice. Right: cold temperature propagating from the injection well in the domain center, showing only temperature values below the threshold of 375 K.
    }
    \label{fig:sequence_A}
\end{figure}

Each model in the sequence is initialized with a reference pressure of 35 MPa (recall that gravitational effects are ignored) and a steady-state temperature of 393 K.
A point source term is attached to the center of each model, with a fixed injection rate of 0.1 $\text{m}^3\text{s}^{-1}$. The injected fluid is assigned the temperature of 313 K. The fixed initial temperature and pressure boundary conditions are prescribed to all the boundaries of the model. The simulated time of each simulation is 27 years. The permeability field and the temperature profile in the middle of one of the simulations are demonstrated in \Cref{fig:sequence_A}.

The following scalar characteristics represent the simulation context, described in \Cref{sec:context}:
\begin{itemize}[noitemsep]
    \item Simulation time step $\Delta t$, [s];
    \item Minimal temperature in the simulation, [K];
    \item Maximum temperature in the simulation, [K];
    \item CFL number: the maximum fluid velocity, multiplied by $\Delta t$ and divided by the characteristic cell length, [-];
    \item Maximum and average values of the enthalpy flux, [J s$^{-1}$];
    \item Maximum and average values of the diffusion flux, [J s$^{-1}$];
    \item Estimate of permeability: the percentage of cells with permeability above the threshold of $10^{-15}$ $\text{m}^2$, [-].
\end{itemize}

The preconditioned linear solver configuration space includes the following decisions:
\begin{itemize}[noitemsep]
    \item GMRES iterative solver is considered. Its \texttt{restart} parameter is tuned within \texttt{\{30, 50, 100\}\texttt};
    \item Three groups of preconditioners are considered: \texttt{\{CPR, System-AMG, Naive\}}. Each group is described below in this list;
    \item The decisions within the \texttt{CPR} preconditioner group:
    \begin{itemize}[noitemsep]
        \item First stage addresses the pressure and temperature submatrices $J_{pp}$ and $J_{TT}$ separately:
        \begin{itemize}[noitemsep]
            \item $J_{pp}$ is approximately solved with a single AMG iteration. Two AMG variations are considered: \texttt{\{BoomerAMG, GAMG\}};
            \item $J_{TT}$ is approximately solved with a single iteration of one of these algorithms: \texttt{\{None, Jacobi, SOR\}}. \texttt{None} corresponds to the classical CPR approach, where $J_{TT}$ is ignored at the first stage;
        \end{itemize}
        \item Second stage addresses the coupled pressure-temperature submatrix with one of these algorithms: \texttt{\{Block-ILU(0), Block-SOR\}};
    \end{itemize}
    \item The decisions within the \texttt{System-AMG} group:
    \begin{itemize}[noitemsep]
        \item The pressure-temperature matrix is approximately solved with a single AMG iteration. Two AMG implementations are considered: \texttt{\{{BoomerAMG}, {GAMG}\}\texttt};
    \end{itemize}
    \item The \texttt{Naive} group incorporates poorly-scalable algorithms. However, they still can compete with the block preconditioners for significantly small or simple models. It includes: \texttt{\{Block-Jacobi, Block-SOR, Block-ILU(0), Block-ILU(1), Block-ILU(2)\}}.
\end{itemize}

The number of available configurations significantly expands when we consider the AMG parameters listed in \Cref{tab:amg_params}. With them, Sequence A considers for solver selection in total 14415 preconditioned linear solver configurations. The vector of solver encoding, which was described in \Cref{sec:solver_space}, contains 31 cells for categorical decisions and 16 cells for numerical decisions.

\subsection{Sequence B: Contact thermo-hydromechanics in fractured porous media with variations in well placements}
This sequence corresponds to the contact-THM problem in 3D porous media with fractures, described in \Cref{sec:thm_model}. The cuboid 3D domain of porous medium is considered, with two clusters of three intersecting fractures. The domain size is $(2 \times 1 \times 1)$ km$^3$.

Each simulation is initialized with a steady-state temperature of 393 K and pressure of 35 MPa. Fixed initial temperature and pressure boundary conditions are prescribed to all the boundaries of the model. The zero displacement boundary condition is prescribed to the bottom of the domain, while the normal compressing force boundary condition is prescribed to the other domain sides: 80 MPa on the top side, 96 MPa on the north and south sides, and 64 MPa on the west and east sides. 
The modeled fluid is water. Unlike Sequence A, the porous medium's porosity and permeability are homogeneous and isotropic, set to $0.05$ and $10^{-13}$ m$^2$, respectively. The shear and bulk moduli are set to $1.2 \cdot 10^{10}$ Pa.

\begin{figure}[t]
    \centering
    \begin{minipage}{0.49\textwidth}
        \centering
        \includegraphics[width=\linewidth]{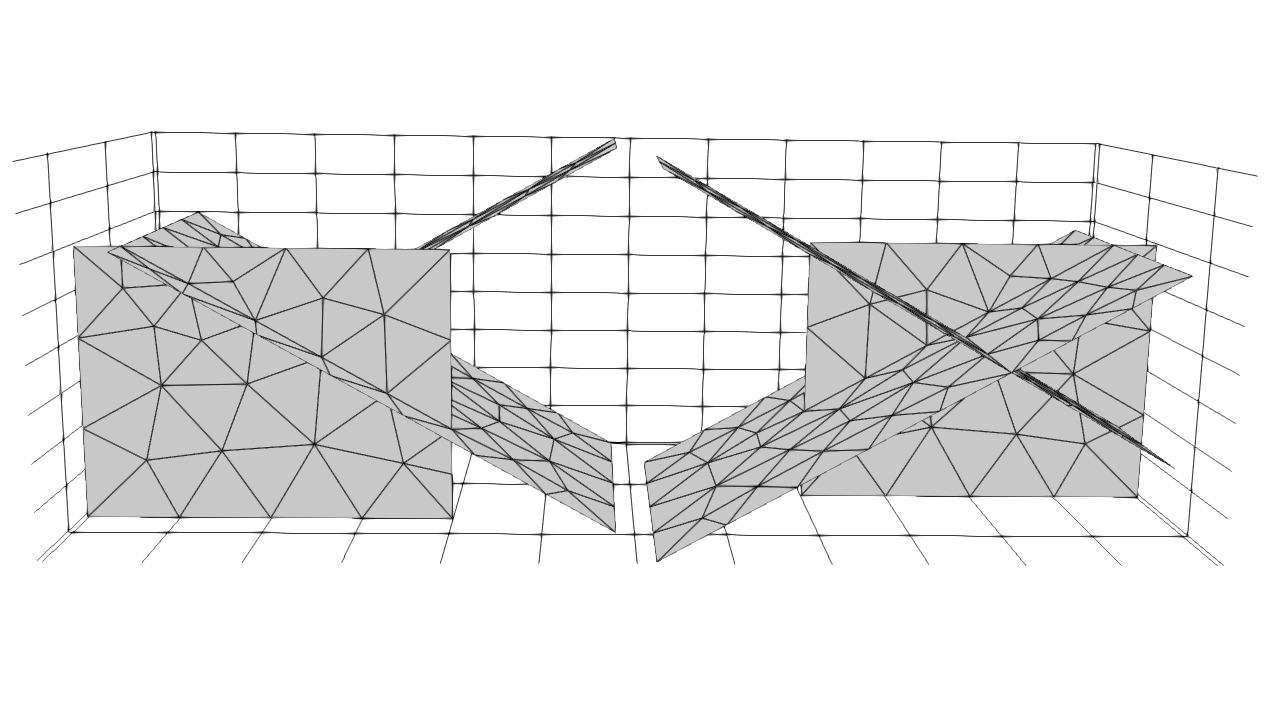}
    \end{minipage}\hfill
    \begin{minipage}{0.49\textwidth}
        \centering
        \includegraphics[width=\linewidth]{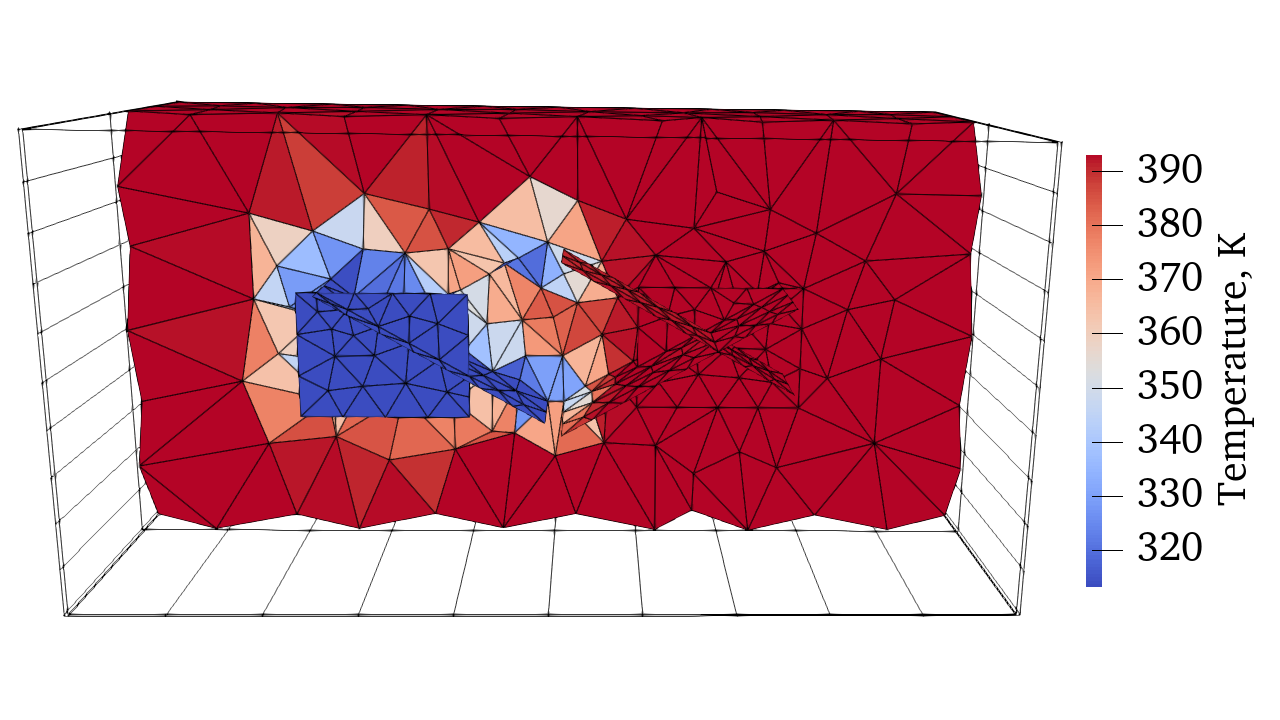}
    \end{minipage}
    \caption{One of the simulations in Sequence B at the end of the simulation time. Left: two clusters of three fractures each. Right: cold temperature propagating from the injection well in the left fracture cluster.
    }
    \label{fig:sequence_B}
\end{figure}

Each simulation features two wells, represented by constant source terms in a fracture cell. One well injects cold fluid of temperature 313 K, and the other produces fluid. The fluid source rate is set to 0.1 $\text{m}^3\text{s}^{-1}$.
The models within this sequence differ from each other by well locations. Five locations are considered for the injection well and five for the production well, resulting in 25 simulations within the sequence. The simulated time of one simulation is 27 years. The fractures and the temperature profile at the end of one of the simulations are demonstrated in \Cref{fig:sequence_B}.

The characteristics that represent the simulation context are similar to \Cref{sec:setup A}, with one difference: the permeability estimate is dropped due to the constant permeability in Sequence B. The preconditioned linear solver configuration space is extended to account for the mechanics subproblem. We also dropped the \texttt{Naive} preconditioners group and only consider \texttt{BoomerAMG} for this problem. Within \texttt{BoomerAMG} parameters, we restricted the ranges of the following parameters to keep the total configurations number reasonable: \texttt{hypre\_boomeramg\_agg\_nl} to \texttt{\{0, 1\}}, \texttt{hypre\_boomeramg\_grid\_sweeps\_all} to \texttt{\{1, 2\}} and 
\texttt{hypre\_boomeramg\_relax\_type\_all} to \texttt{\{SOR/Jacobi, Jacobi\}}. The description of the configuration space is given by:
\begin{itemize}[noitemsep]
    \item GMRES iterative solver is considered. This time, \texttt{restart} is fixed to 100;
    \item The mechanics subsolver is based on \texttt{BoomerAMG}. The fixed-stress approximation is applied to decouple mechanics;
    \item Two subsolver groups are considered for the flow and transport subproblem: \texttt{\{CPR, System-AMG\}}. Each group is described below in this list;
    \item The decisions within the \texttt{CPR} subsolver group:
    \begin{itemize}[noitemsep]
        \item First stage addresses the pressure and temperature submatrices $J_{pp}$ and $J_{TT}$ separately:
        \begin{itemize}[noitemsep]
            \item $J_{pp}$ is approximately solved with a single AMG iteration. Only \texttt{BoomerAMG} implementation is considered;
            \item $J_{TT}$ is approximately solved with a single iteration of one of these algorithms: \texttt{\{None, SOR\}}. \texttt{None} corresponds to the classical CPR approach, where $J_{TT}$ is ignored at the first stage;
        \end{itemize}
        \item Second stage addresses the coupled pressure-temperature submatrix with one of these algorithms: \texttt{\{Block-ILU(0), Block-SOR\}};
    \end{itemize}
    \item The decisions within the \texttt{System-AMG} subsolver group:
    \begin{itemize}[noitemsep]
        \item The pressure-temperature matrix is approximately solved with a single AMG iteration. Only \texttt{BoomerAMG} implementation is considered.
    \end{itemize}
\end{itemize}
In total, Sequence B considers 32000 preconditioned linear solver configurations for solver selection. The vector of solver encoding described in \Cref{sec:solver_space} contains 21 cells for categorical decisions and 10 cells for numerical decisions.

\section{Numerical Experiments}
\label{sec:experiments}

In the numerical experiments, we run simulation sequences A and B described in \Cref{sec:target_simulations} and evaluate the solver selection algorithm: how much it improves simulation performance over a baseline configuration, how robust and reproducible the improvement is, and how much the size of the performance dataset impacts this improvement. The experiments are organized into three steps.

\Cref{sec:collecting_statistics} collects representative statistics about the simulated sequences and how the solver configurations perform with them. The goal is to estimate the problem complexity and understand how much improvement we should expect from solver selection. Note that this step is not meant to be done when applying the solver selection algorithm in practice.

\Cref{sec:solver_selection_experiment} runs the simulation sequences using the solver selection algorithm in realistic conditions: performance data needs to be collected during the run time. Within this experiment, we explore how the selection policy improves as we collect more performance data. We also illustrate how the simulation state affects the expected run time and demonstrate how the collected performance data can be analyzed to refine the solver configuration ranges for future simulations. Next, we present the machine learning run time overhead and analyze the run time results achieved by the solver selection algorithm.

In \Cref{sec:experiment_against_optimal_solver}, we re-run the simulation sequences, this time using all the previously collected performance data. This approach aims to achieve optimal simulation performance but requires a large performance dataset, making it infeasible in practical simulations. We compare its performance with the practical case from \Cref{sec:solver_selection_experiment}.

The solver selection algorithm is non-deterministic due to the random selections within the initial exploration and random noise in time measurements. Therefore, we repeat each sequence run five times and report average results and standard deviations. To further test robustness, we shuffle the setups (the slicing in Sequence A, the well placements in Sequence B) within each run. 
Different random seeds are used for shuffling, selection, and the machine-learning models in each of the five runs.

All simulations are run on a single core.
We use the SciKit-Learn implementation of the gradient boosting classification and regression algorithms \cite{scikit-learn}. 

\subsection{Collecting Statistics}
\label{sec:collecting_statistics}
To collect statistics, we run the simulations in Sequences A and B with random solver configurations, similarly to the initial exploration described in \Cref{sec:initial_exploration}, but throughout the whole sequence. The solver configuration is switched randomly after each linear system, regardless of whether it was solved successfully or not. This process is repeated five times with the simulations in shuffled order.

\begin{figure}[ht]
\centering
\begin{minipage}{0.49\textwidth}
\begin{table}[H]
\centering
\setlength{\abovecaptionskip}{5pt}  %
\begin{tabular}{lr}
\multicolumn{2}{c}{\textbf{Sequence A Random}} \\
\toprule
Num. solver configurations & 16575 \\
Num. data points & 19573 \\
Configurations tried, \%  & 69.2 \\
Success rate, \% & 64.5 \\
Always success, \%  & 63.2 \\
Always failure, \%  & 35.0 \\
\hdashline
Run time average, s & 3.70 \\
Run time median, s & 2.54 \\
Run time min, s & 1.17 \\
Run time max, s & 178.30 \\
\bottomrule
\end{tabular}
\end{table}
\end{minipage}
\hfill
\begin{minipage}{0.49\textwidth}
\vspace{6mm}
\begin{figure}[H] %
\centering
\includegraphics{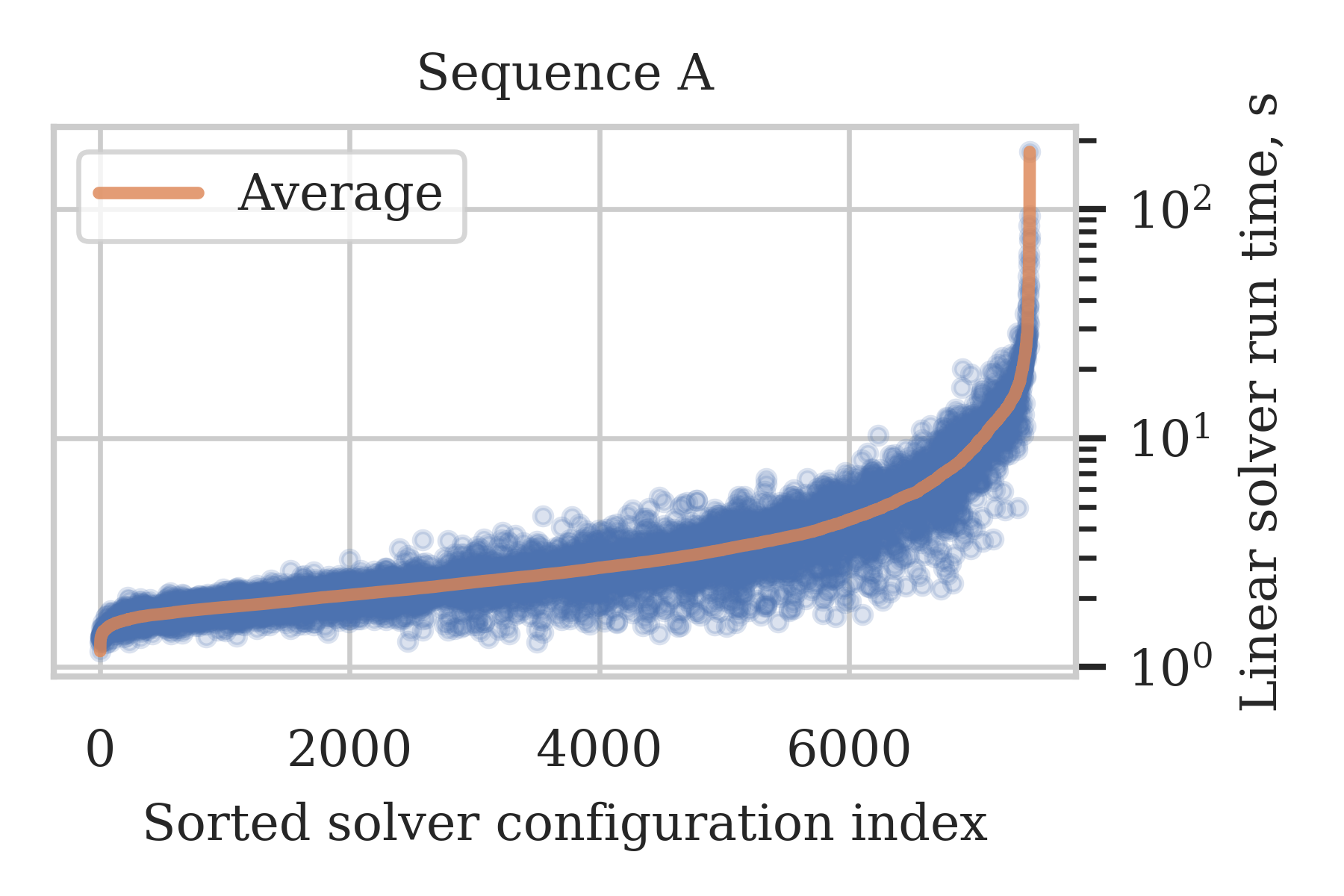}
\end{figure}
\end{minipage}
\begin{minipage}{0.49\textwidth}
\begin{table}[H]
\centering
\setlength{\abovecaptionskip}{5pt}  %
\begin{tabular}{lr}
\multicolumn{2}{c}{\textbf{Sequence B Random}} \\
\toprule
Num. solver configurations & 32000 \\
Num. data points & 21661 \\
Configurations tried, \% & 49.4 \\
Success rate, \% & 91.3 \\
Always success, \%  & 90.5 \\
Always failure, \%  & 8.0 \\
\hdashline
Run time average, s & 5.96 \\
Run time median, s & 5.47 \\
Run time min, s & 3.56 \\
Run time max, s & 18.35 \\
\bottomrule
\end{tabular}
\end{table}
\end{minipage}
\hfill
\begin{minipage}{0.49\textwidth}
\vspace{6mm}
\begin{figure}[H] %
\centering
\includegraphics{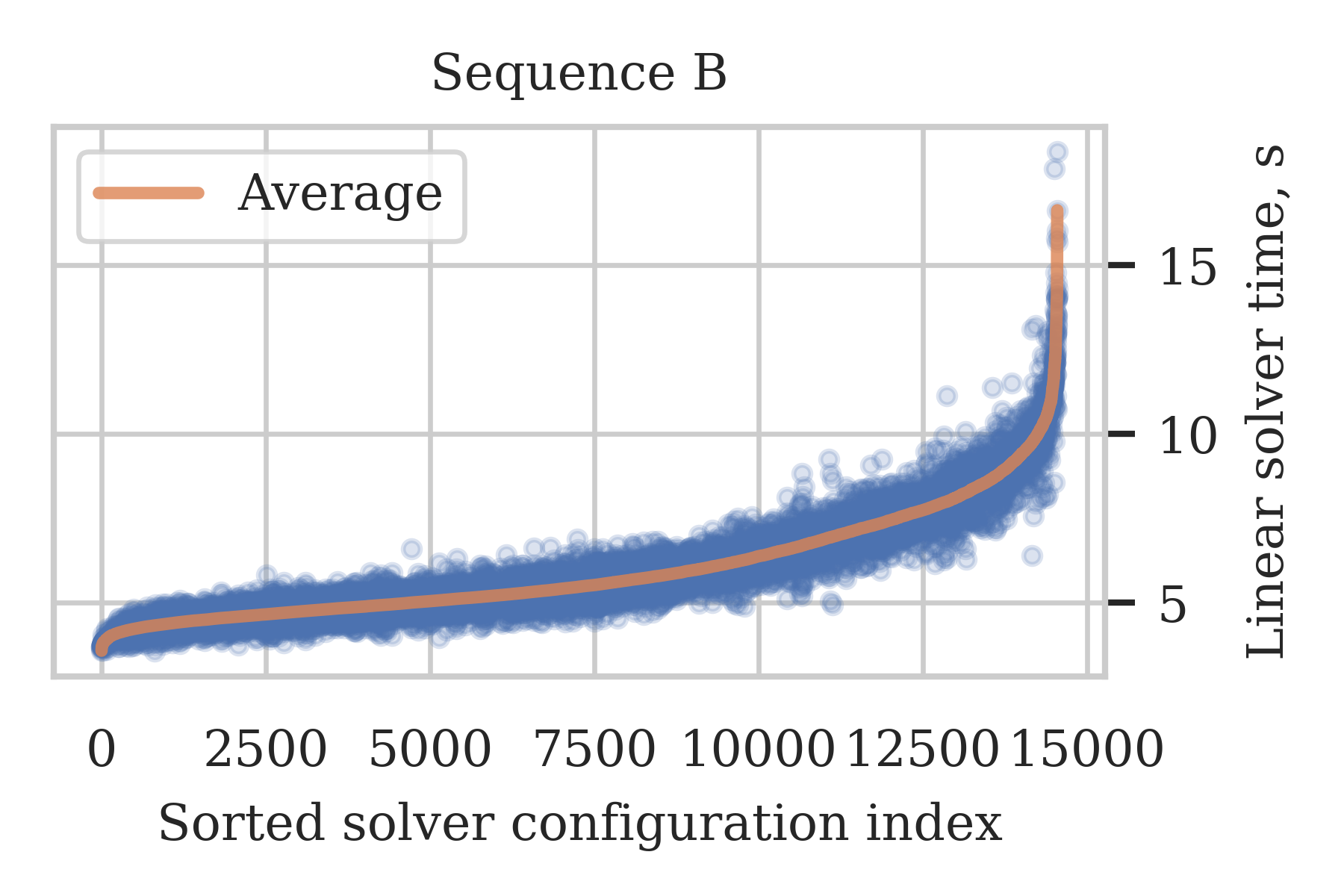}
\end{figure}
\end{minipage}
\caption{Performance statistics collected during the experiments with the random solver configurations for Sequences A and B. ``Always success'' and ``Always failure'' indicate the percentage of solver configurations that consistently solved or failed to solve the given linear systems (for clarity, some configurations may not have been applied more than once). The percentages of ``Always success'' and ``Always failure'' are calculated based on the total number of attempted configurations. Their sum does not equal 100\% due to a few inconsistent solvers. The run time data corresponds to the time it took to construct the linear solver and solve the given problem. The run time statistics only consider the successful linear solve attempts. The data points in the right-hand figures represent all successful linear solve attempts, with solver configurations indexed by their sorted average run time.}
\label{fig:statistics_random_runs}
\end{figure}

The performance statistics are shown in \Cref{fig:statistics_random_runs}, with summary statistics in the tables and run time distributions in the figures, where we show the solver configurations' performance separately from each other. For this, we indexed the solver configurations and sorted them based on average performance. This representation completely ignores the varying simulation states (and, consequently, the varying complexity of the linear systems). We return to this point later.

Due to the large number of possible options, not every solver configuration was tried in the simulations. However, the sampling was random, so the statistics are representative. The major part of the solver configurations solved the given problems successfully. Sequence A contains more bad configurations than Sequence B, including a larger number that always fail and several that converge but require orders of magnitude more time than that of the average solver.
We note that Sequence A includes only 15 clearly bad solver configurations from the \texttt{Naive} group. The remaining failing configurations correspond to various suboptimal System-AMG or CPR setups, which were not obviously bad a priori.

\subsection{Solver Selection Experiment}
\label{sec:solver_selection_experiment}

We run Sequences A and B with the solver selection algorithm. Each sequence run begins with no performance data. The initial exploration phase described in \Cref{sec:initial_exploration} is limited by \texttt{num\_initial = 64} linear systems, and the maximum number of allowed failures for one linear system is 20. The batch size for updating the machine learning pipeline described in \Cref{sec:incremental_learning} is set to 64. The solver selection algorithm is allowed to select a configuration for each Newton iteration. The whole experiment is repeated five times with the simulations in shuffled order.

\begin{figure}[htb]
    \centering
    \begin{minipage}{0.49\textwidth}
        \centering
        \includegraphics[width=\linewidth]{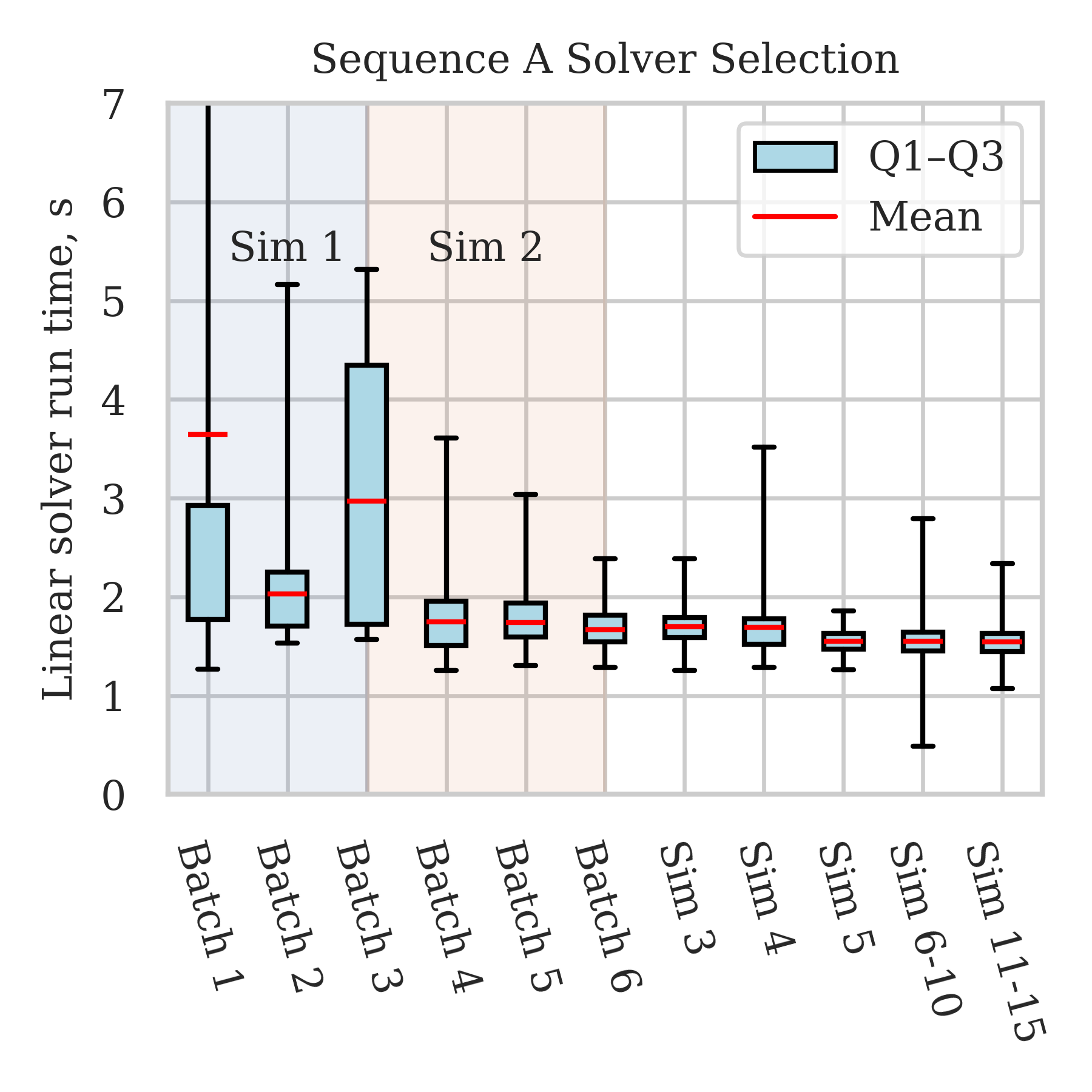}
    \end{minipage}\hfill
    \begin{minipage}{0.49\textwidth}
        \centering
        \includegraphics[width=\linewidth]{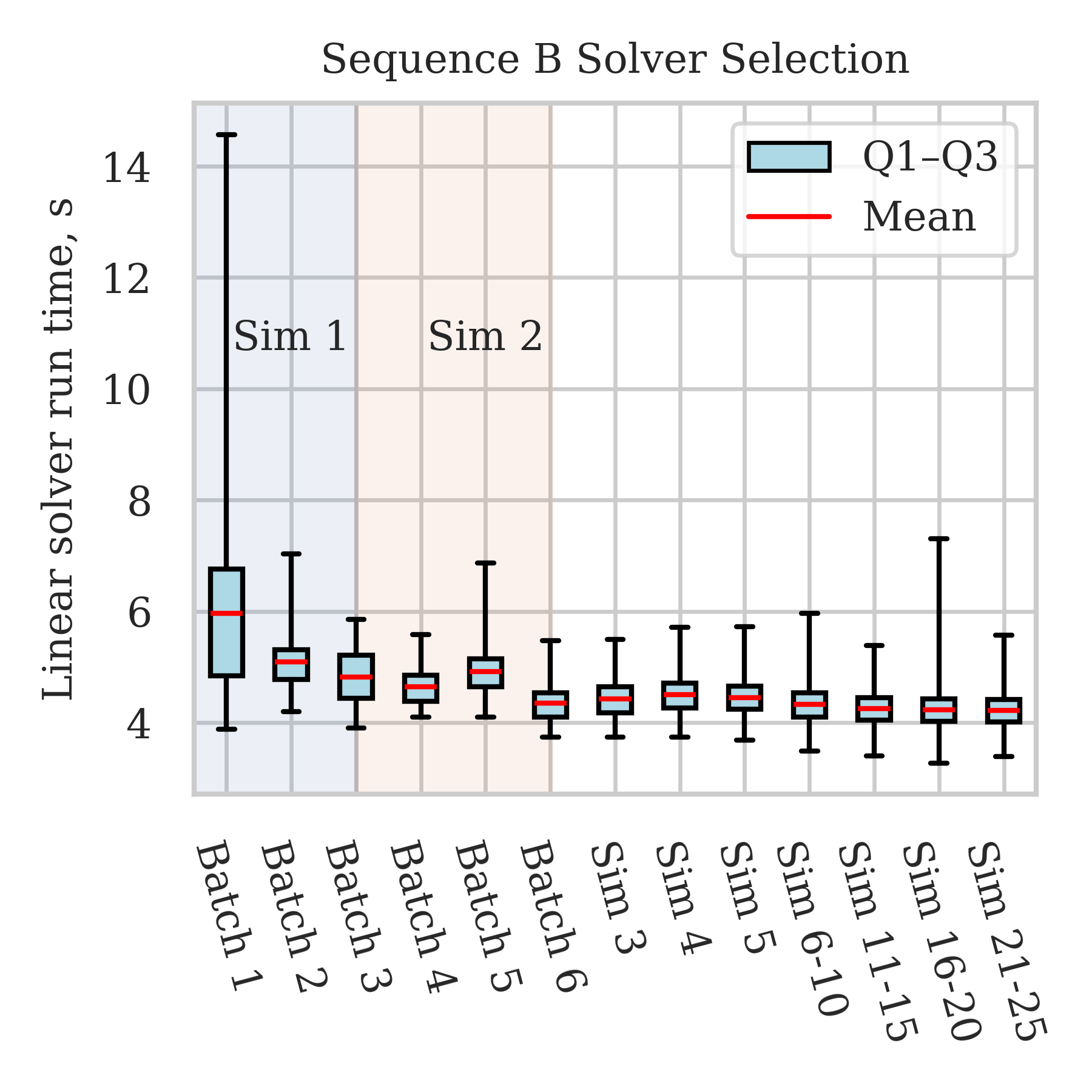}
    \end{minipage}
    \caption{Distributions of the linear solver run times in successive simulations of Sequences A and B. 
    The first two simulations are described in detail. Datasets are grouped by machine learning batches. Batch 1 covers the initial exploration. Batch 3 spans from Simulation 1 to Simulation 2.
    From Simulation 3 onward, distributions are shown per simulation, and from Simulation 6 they are also aggregated by groups of simulations.
    Only the successful solution attempts are included. The red lines denote the mean value, the blue boxes denote the interquartile range, spanning from the 25th percentile (Q1) to the 75th percentile (Q3). The black vertical lines show the minimum and maximum values of each distribution. The figure for Sequence A is truncated on the top, as the maximum run time in Batch 1 is 75 seconds.}
    \label{fig:boxplot_sim_number_batches}
\end{figure}

Simulation performance is shown in \Cref{fig:boxplot_sim_number_batches}, focusing on the early part of the sequence, as performance stabilizes later. The first two simulations of each sequence are shown in detail, divided into three sections corresponding to the machine-learning pipeline’s batch updates. The first batch corresponds to the initial exploration, where decisions are taken randomly. Batches 2 and 3 cover the remainder of the first simulation, with batch 3 extending into the second simulation, and batches 3 -- 6 fully cover simulation 2. 
After the initial exploration, the spread within Q1-Q3 in Sequence A remains wide until simulation 4, after which the mean stabilizes and the spread narrows. Occasional outliers are still present until the end of the sequence, but their effect on the average performance is negligible. In Sequence B, the spread within Q1-Q3 decreases and the behavior stabilizes starting at simulation 3. After that, only a few outliers are present, showing that the solver selection algorithm consistently avoids bad decisions.

\begin{figure}[htb]
    \centering
    \begin{minipage}{0.49\textwidth}
        \centering
        \includegraphics[width=\linewidth]{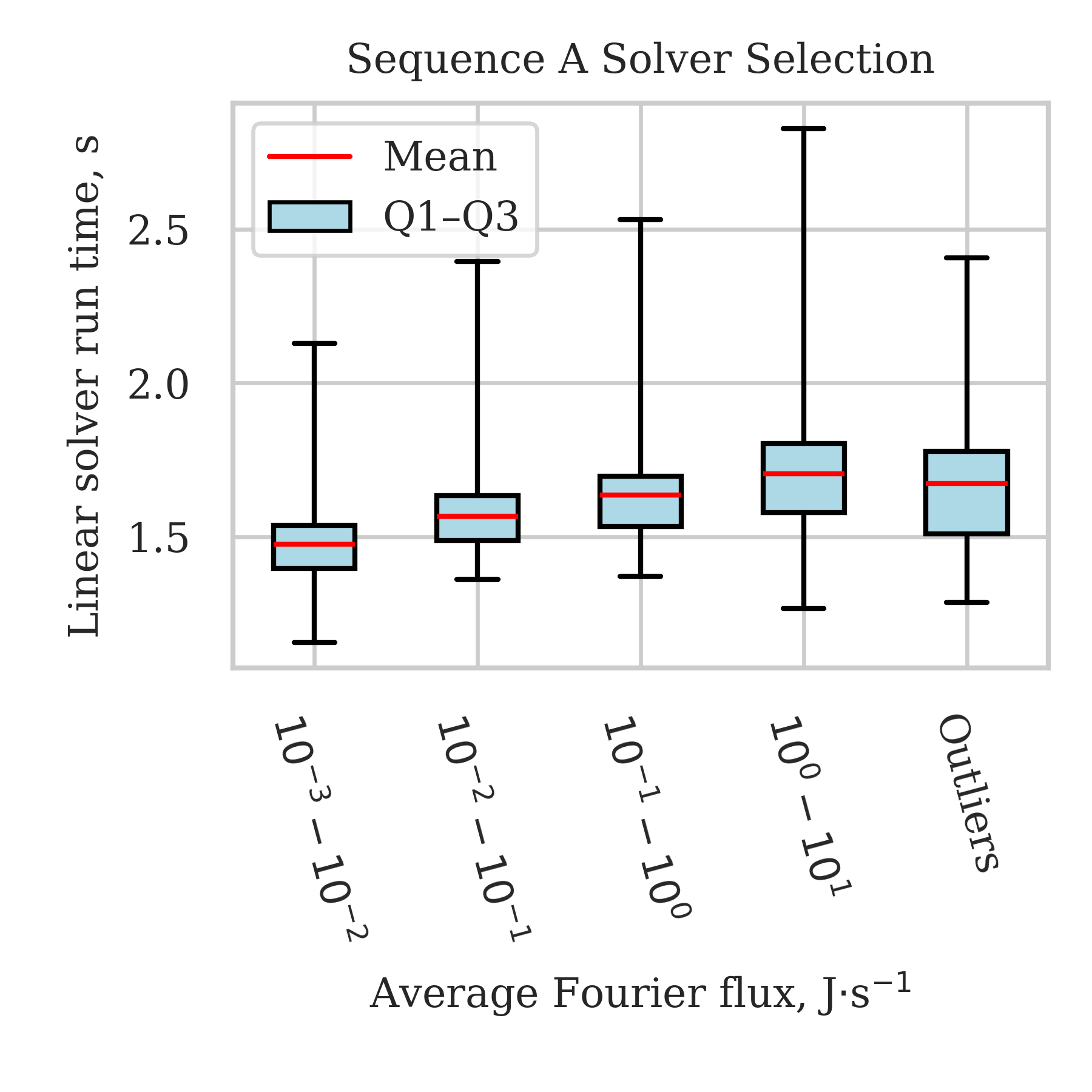}
    \end{minipage}\hfill
    \begin{minipage}{0.49\textwidth}
        \centering
        \includegraphics[width=\linewidth]{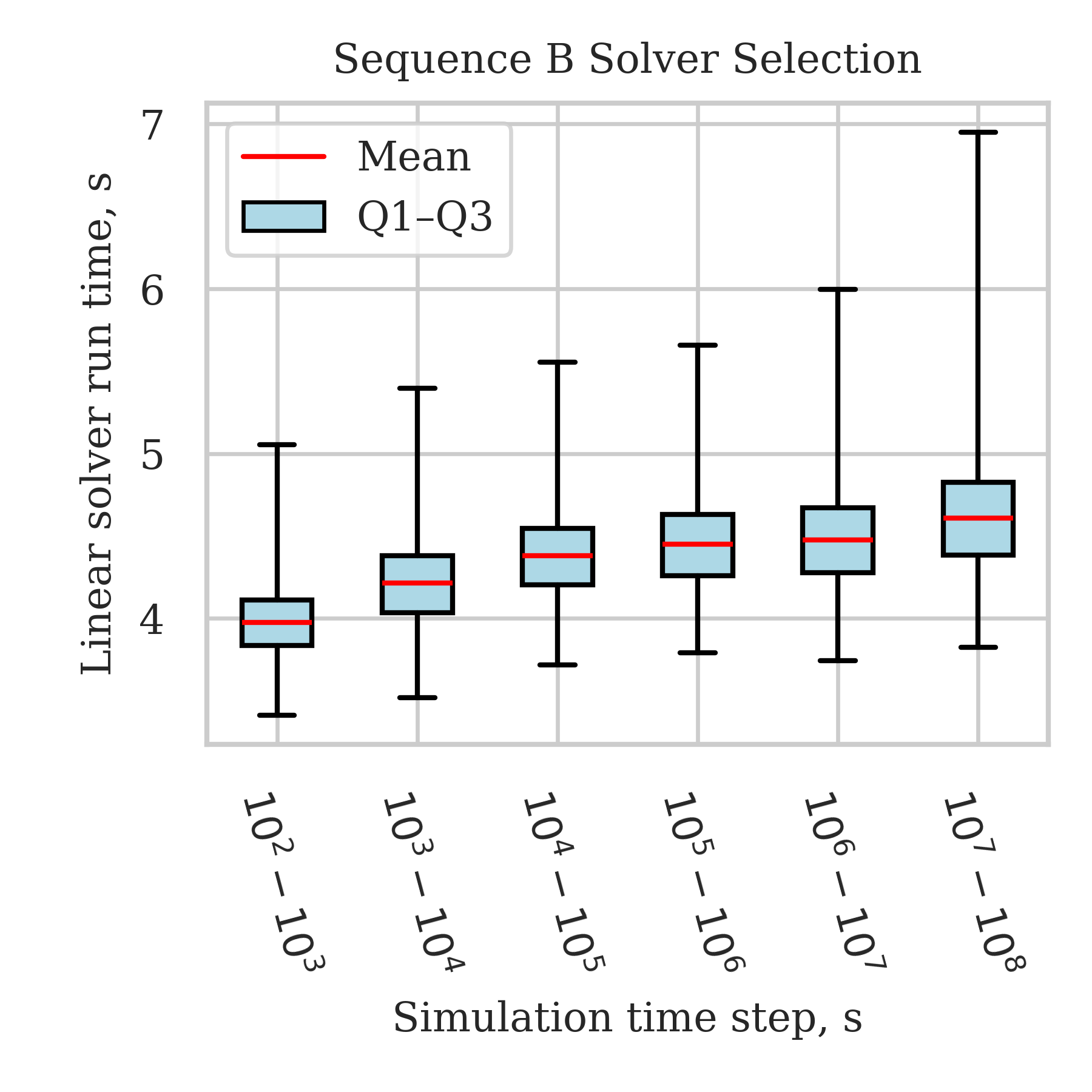}
    \end{minipage}
    \caption{Distributions of the linear solver run times subject to simulation context. The datasets correspond to simulations 11 -- 15 and 16 -- 25 of Sequences A and B, respectively. 
    Only the successful solution attempts are included. The red lines denote the mean value, the blue boxes denote the interquartile range, spanning from the 25th percentile (Q1) to the 75th percentile (Q3). The black vertical lines show the minimum and maximum values of each distribution.
    }
    \label{fig:boxplot_context}
\end{figure}

We analyze the simulation performance as a function of the simulation context in \Cref{fig:boxplot_context}. The boxes based on the average Fourier flux and the simulation time step are used to group the data of Sequences A and B, respectively. The Fourier flux values contain about 0.8\% outliers, all of which fall outside the plotted boxes, so we place them in a separate box. The figure shows that simulation context affects the performance, and is thus necessary for accurate performance predictions. This also illustrates why the initial exploration can be insufficient: the machine learning algorithm encounters previously unseen simulation contexts in batches 2 and 3 of \Cref{fig:boxplot_sim_number_batches}, leading to a large spread and a suboptimal average performance. However, the later batches demonstrate increased accuracy as more simulation states are encountered.

\begin{figure}[htb]
    \centering
    \begin{minipage}{0.49\textwidth}
        \centering
        \includegraphics[width=\linewidth]{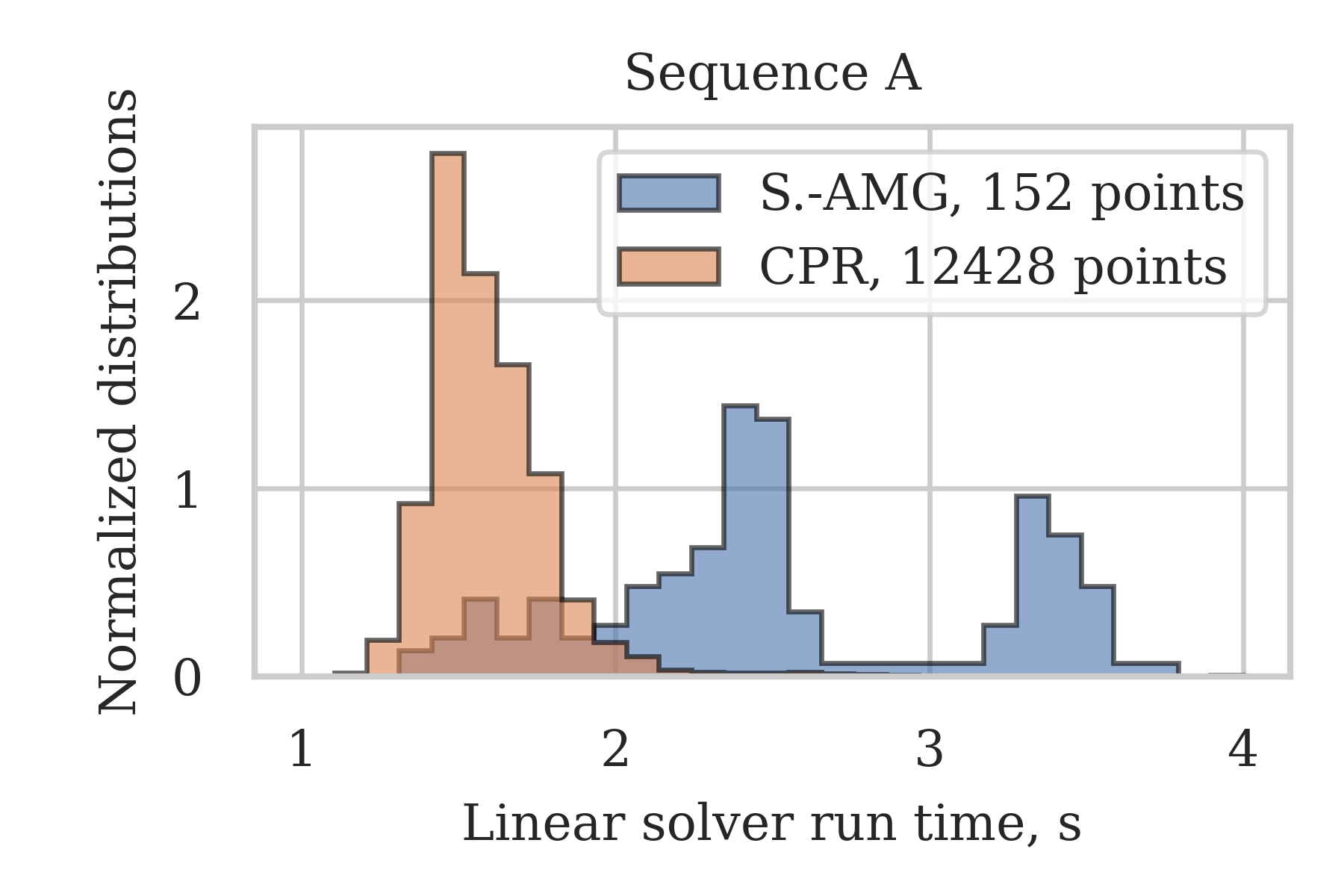}
    \end{minipage}\hfill
    \begin{minipage}{0.49\textwidth}
        \centering
        \includegraphics[width=\linewidth]{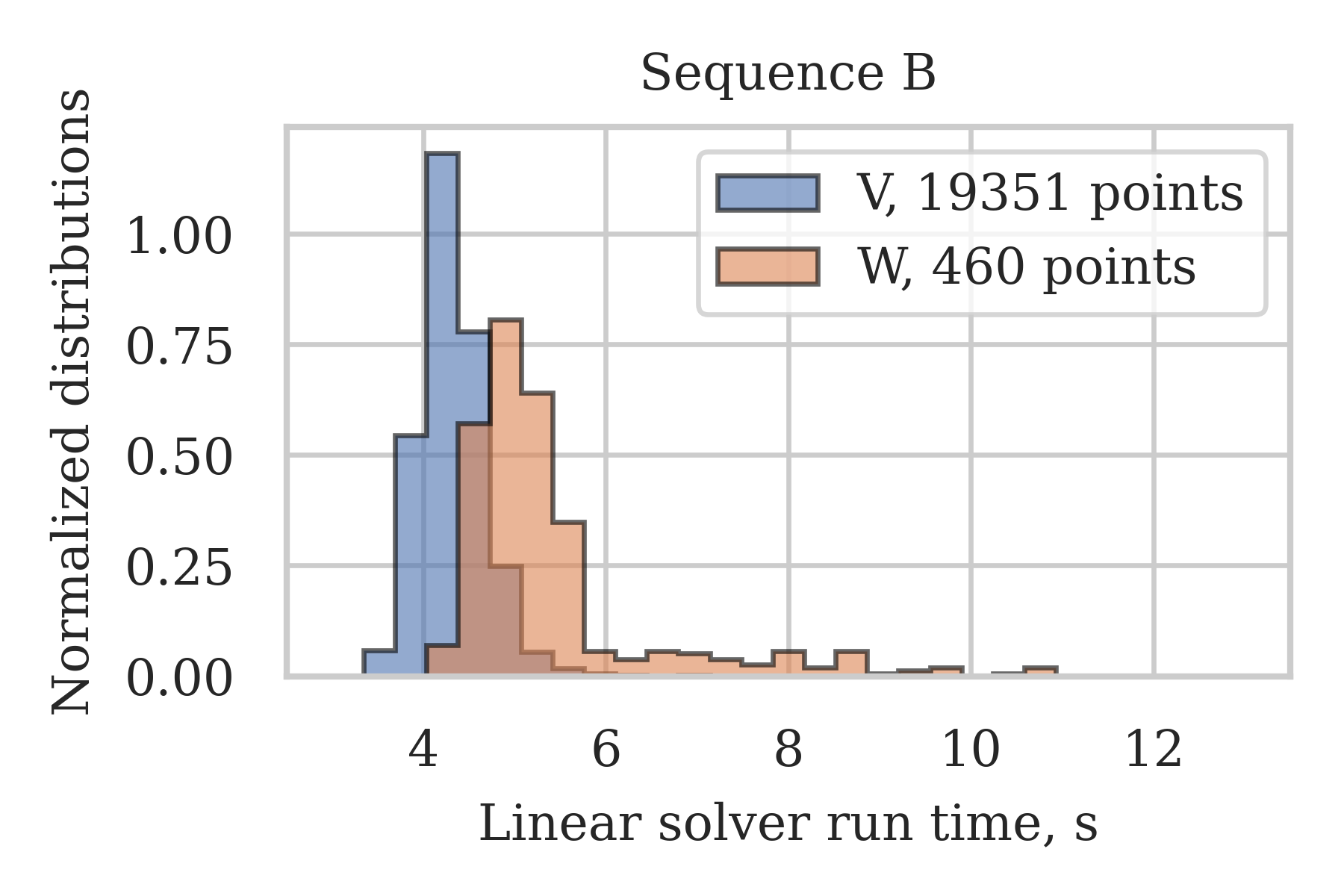}
    \end{minipage}
    \caption{Distributions of the linear solver run times resulting from particular algorithmic choices: System-AMG or CPR for Sequence A and V- or W-cycle of the mechanics AMG subsolver for Sequence B.    
    Run time ranges are truncated to remove outliers. The legend also specifies the number of data points in the distribution.
    }
    \label{fig:analyzing_decisions}
\end{figure}

\Cref{fig:analyzing_decisions} shows how algorithmic choices impact the linear solver run time. Namely, we demonstrate the choices: (i) System-AMG or CPR for Sequence A; (ii) AMG V- or W-cycle for the mechanics subsolver of Sequence B. Only the successful solution attempts are included. 
As the number of data points differs between the choices (e.g., V- or W-cycle), the data is shown as normalized histograms, so the total area of each distribution equals 1.
Performance may depend on complex, nonlinear combinations of decisions, requiring more advanced data analysis and visualization methods. However, the chosen features clearly separate the distributions and serve as a good illustration of how a similar analysis, conducted after the simulations, allows us to interpret the solver selection algorithm decisions and to refine the configuration space for future simulations.

\begin{figure}[htb]
    \centering
    \begin{minipage}{0.49\textwidth}
        \centering
        \includegraphics[width=\linewidth]{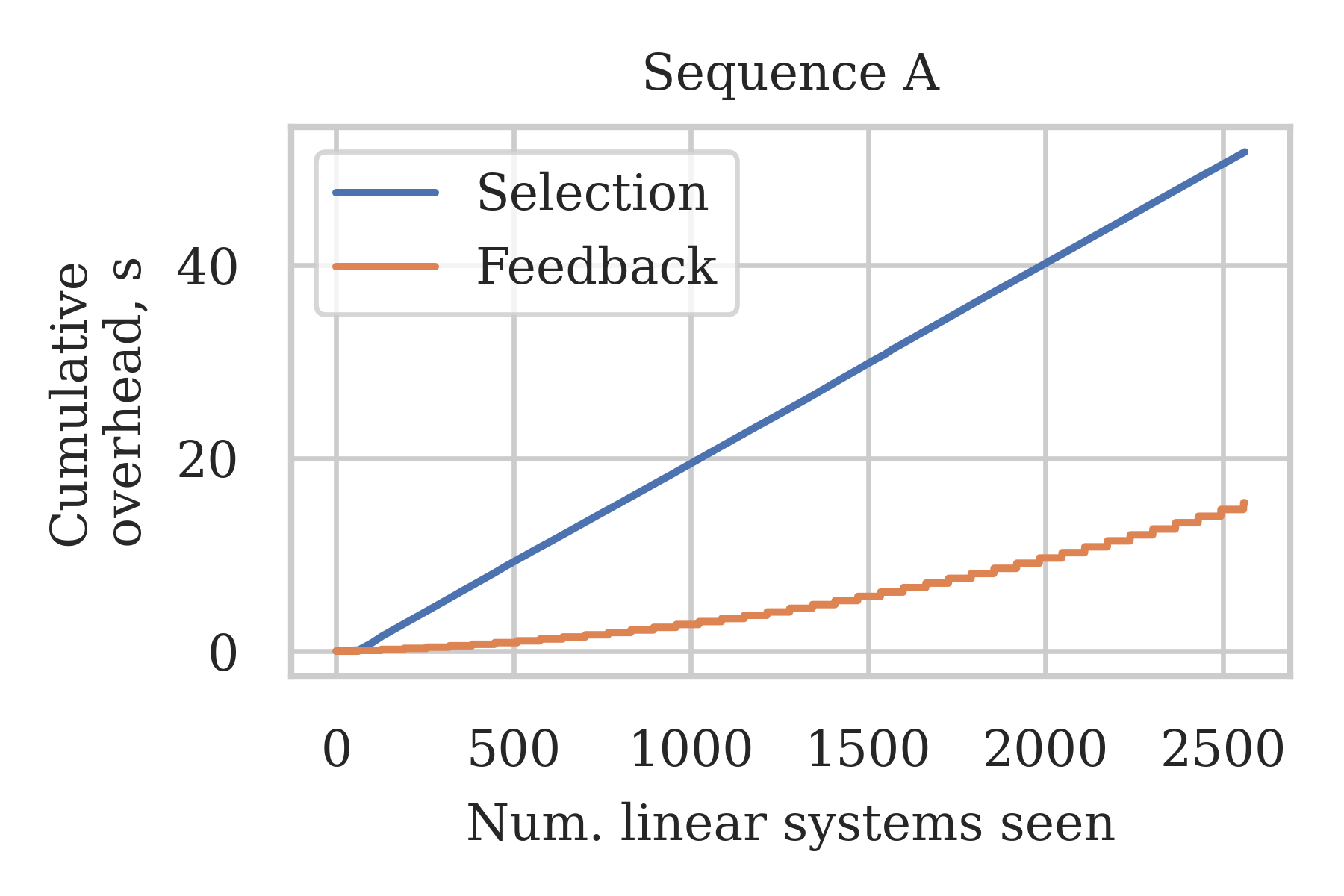}
    \end{minipage}\hfill
    \begin{minipage}{0.49\textwidth}
        \centering
        \includegraphics[width=\linewidth]{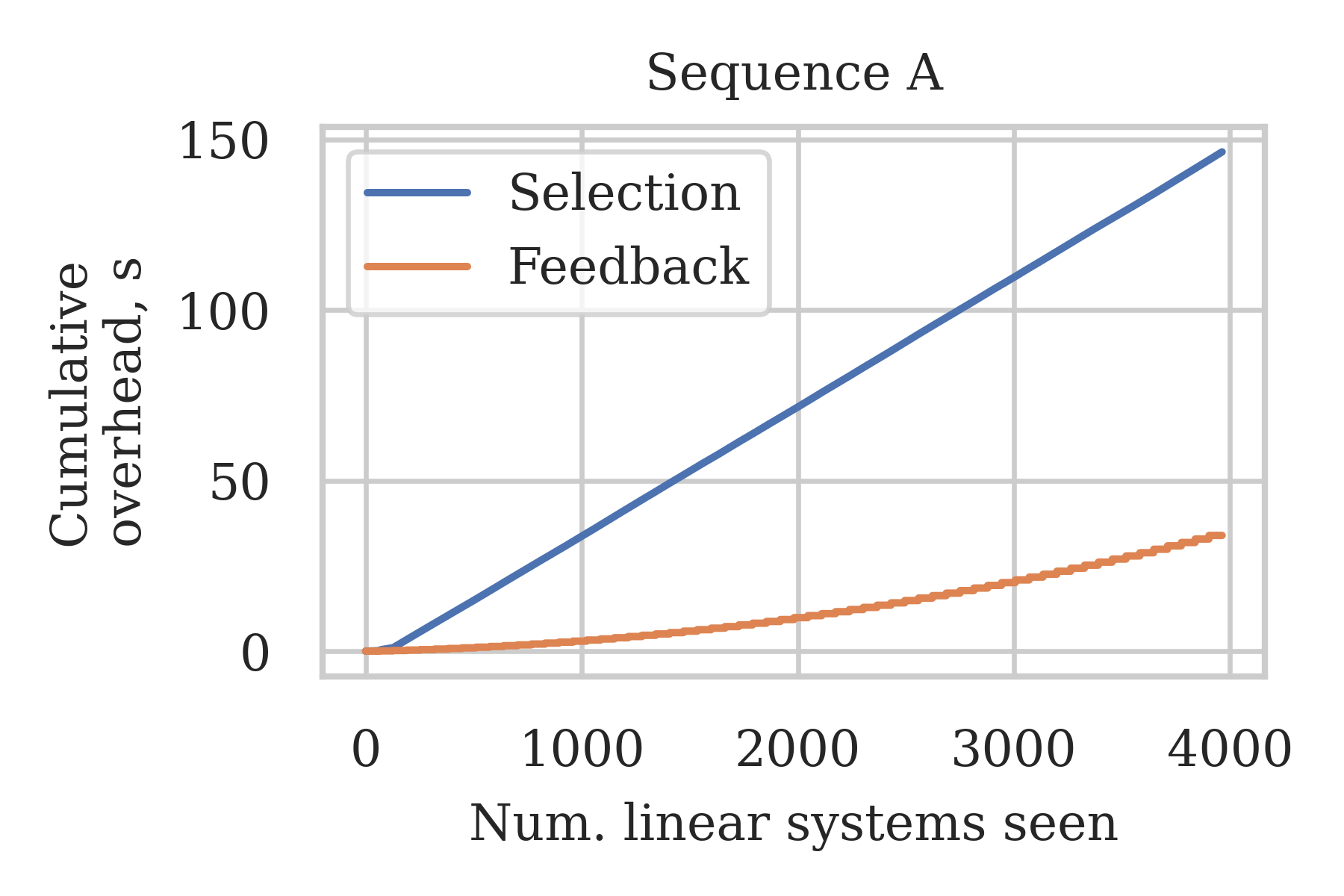}
    \end{minipage}
    \caption{Computation overhead caused by the solver selection algorithm in Sequences A and B.}
    \label{fig:ml_overhead}
\end{figure}

The machine learning pipeline within the solver selection algorithm inevitably introduces computational overhead, as shown in \Cref{fig:ml_overhead}. The major part occurs during the selection, which involves the machine learning pipeline prediction for 16 and 32 thousands solver configurations for every new linear system. On average, it takes 0.020 ± 0.003 s and 0.037 ± 0.006 s per linear system for Sequences A and B, respectively, and this value does not increase with the number of data points. The sawtooth growing pattern in feedback overhead is attributed to the batched updates. We observe near-linear scaling of time to refit the gradient boosting model (result not provided here). As we conclude, the overhead is negligible relative to the time required to solve linear systems.

\begin{figure}[hbt]
\centering
\begin{minipage}{0.49\textwidth}
\begin{table}[H]
\centering
\setlength{\abovecaptionskip}{5pt}  %
\begin{tabular}{lr}
\multicolumn{2}{c}{\textbf{Sequence A Solver Selection}} \\
\toprule
Num. solver configurations & 16575 \\
Num. data points & 14169 \\
Configurations tried, \%  & 3.40 \\
Success rate after init. expl., \% & 88.88 \\
\hdashline
Run time average, s & 1.65 \\
Run time median, s & 1.56 \\
Run time min, s & 0.48 \\
Run time max, s & 75.27 \\
\bottomrule
\end{tabular}
\end{table}
\end{minipage}
\hfill
\begin{minipage}{0.49\textwidth}
\vspace{6mm}
\begin{figure}[H] %
\centering
\includegraphics{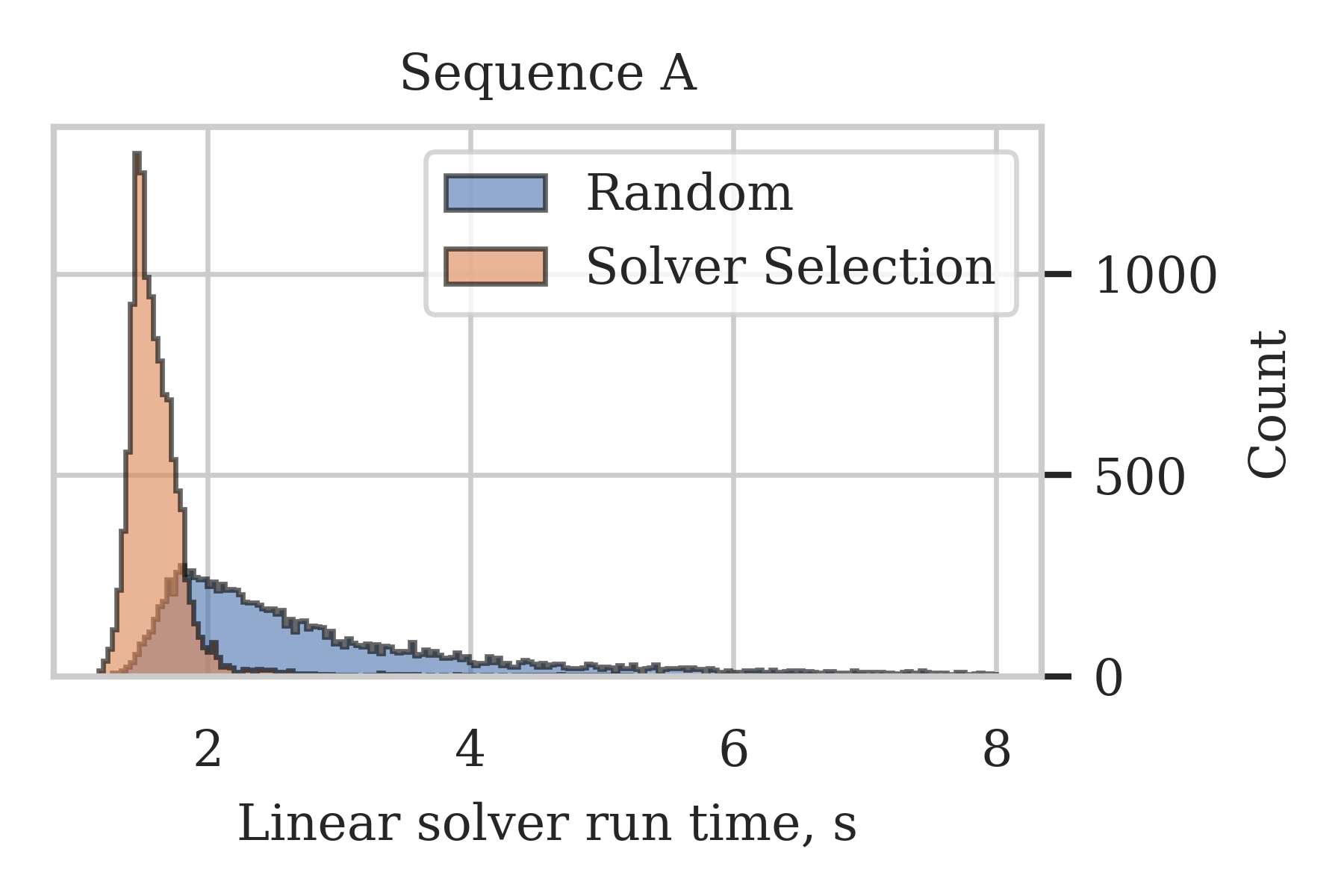}
\end{figure}
\end{minipage}
\begin{minipage}{0.49\textwidth}
\begin{table}[H]
\centering
\setlength{\abovecaptionskip}{5pt}  %
\begin{tabular}{lr}
\multicolumn{2}{c}{\textbf{Sequence B Solver Selection}} \\
\toprule
Num. solver configurations & 32000 \\
Num. data points & 19830 \\
Configurations tried, \% & 1.88 \\
Success rate after init. expl., \% & 99.92 \\
\hdashline
Run time average, s & 4.34 \\
Run time median, s & 4.29 \\
Run time min, s & 3.27 \\
Run time max, s & 14.57 \\
\bottomrule
\end{tabular}
\end{table}
\end{minipage}
\hfill
\begin{minipage}{0.49\textwidth}
\vspace{6mm}
\begin{figure}[H] %
\centering
\includegraphics{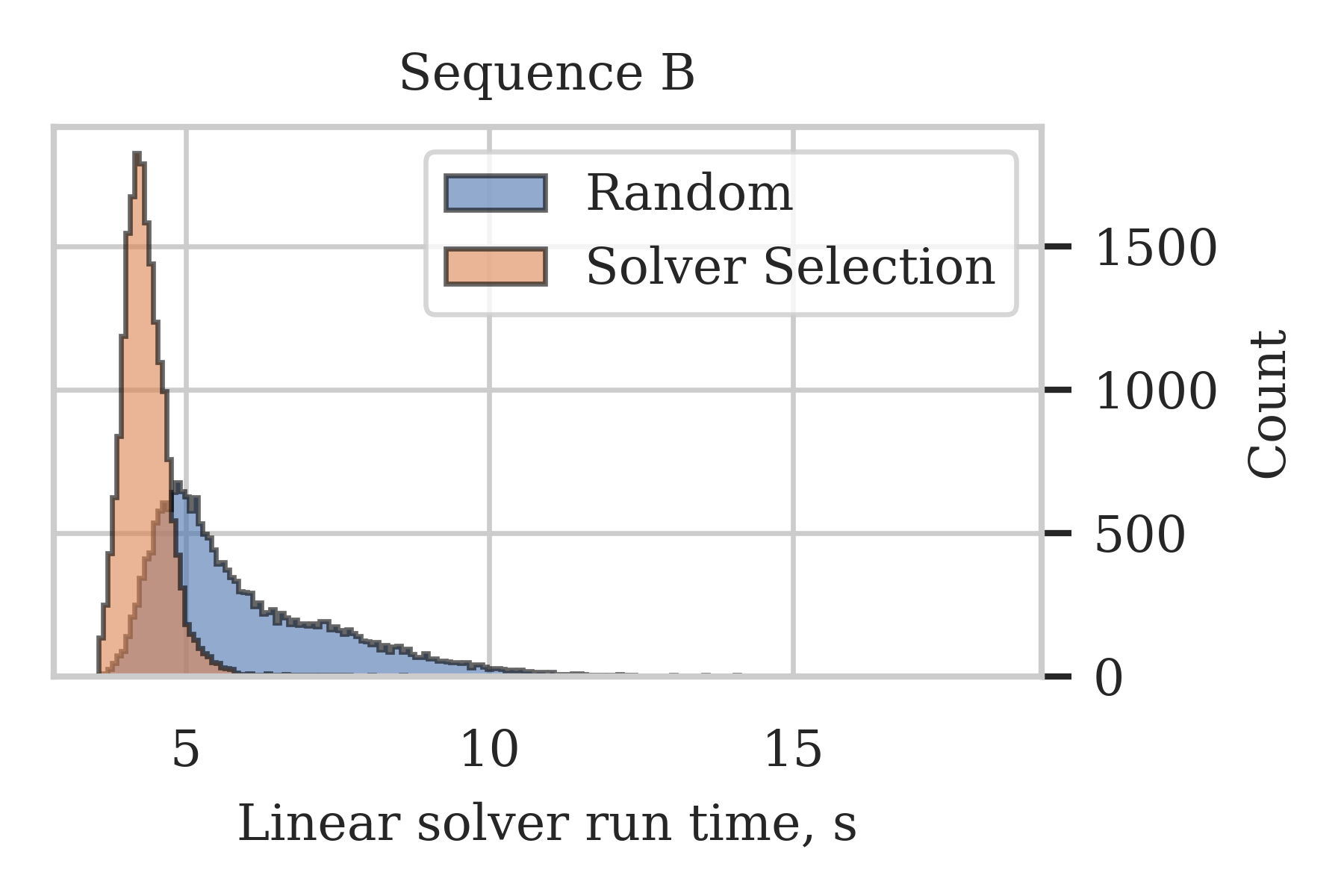}
\end{figure}
\end{minipage}
\caption{Performance statistics of the experiment with the solver selection algorithm for Sequences A and B, including the initial exploration data. The left-hand side tables are similar to those in \Cref{fig:statistics_random_runs}.
The success rate is computed for the solver selection algorithm decisions, excluding initial random exploration.
The right-hand side histograms compare the simulation performance in this experiment with the data from \Cref{fig:statistics_random_runs}. The histogram for Sequence A is truncated on the right. The random experiment dataset contains 939 points with the run time > 8 seconds, whereas the solver selection dataset contains only 19 such points. The histogram for Sequence B is not truncated.
}
\label{fig:statistics_solver_selection}
\end{figure}

\Cref{fig:statistics_solver_selection} summarizes the experiment results and compares them to the statistics collected in \Cref{sec:collecting_statistics}. For both sequences, fewer than 2\% of all available solver configurations were used. However, the solver selection algorithm was able to determine well-performing configurations among them and stick to them until the end of the simulations, increasing the success rate above 88\% and 99\% for Sequences A and B, respectively. The histograms show that the solver selection algorithm can robustly avoid bad decisions and stay within the left side of the random experiment distribution's bell curve.
The mean successful linear solver run time in this experiment corresponds to the configurations within the top 1.9\% and 2.7\% for Sequences A and B, respectively, in the histograms of \Cref{fig:statistics_random_runs}. 
The cost of this result is solving only 64 linear systems randomly during the initial exploration, along with negligible machine learning overhead.

\subsection{Comparing Against Optimal Solver}
\label{sec:experiment_against_optimal_solver}

\begin{figure}[hbt]
    \centering
    \begin{minipage}{0.49\textwidth}
        \centering
        \includegraphics[width=\linewidth]{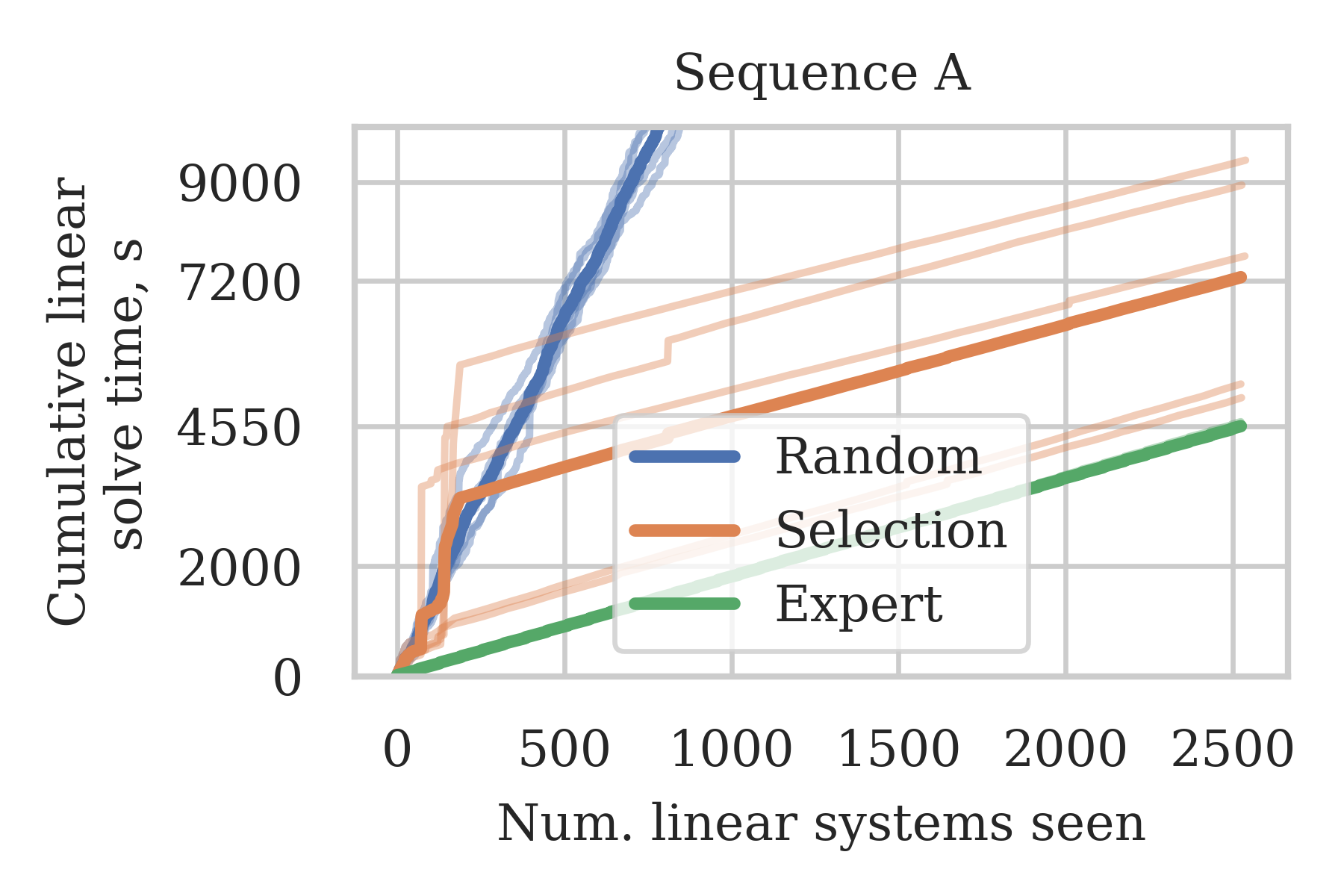}
    \end{minipage}\hfill
    \begin{minipage}{0.49\textwidth}
        \centering
        \includegraphics[width=\linewidth]{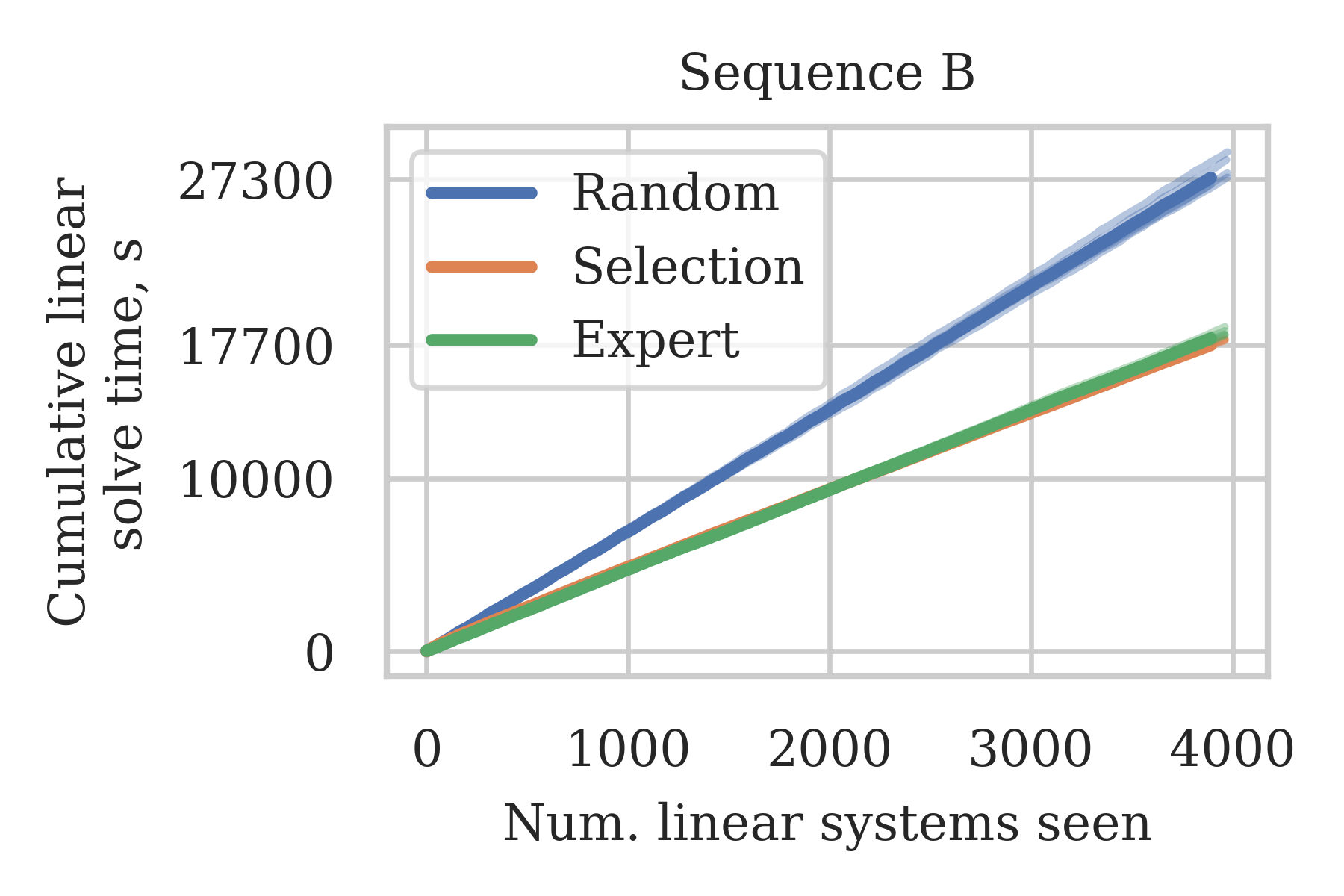}
    \end{minipage}
    \caption{Cumulative run time of linear solvers in Sequences A and B. ``Random'', ``Selection'' and ``Expert'' correspond to the runs described in \Cref{sec:collecting_statistics}, \Cref{sec:solver_selection_experiment} and \Cref{sec:experiment_against_optimal_solver}, respectively. Five faint lines of the same color correspond to five restarts of each experiment, with the shuffled order of the simulations and different random seeds. Bold lines show the mean across five runs. The ``Random'' case in Sequence A is truncated on the top.
    }
    \label{fig:runtime_expert}
\end{figure}

In the previous experiments, we have collected numerous data points, enough to determine the optimal solver configurations for the considered simulation setups. We can plug in all this data to the solver selection algorithm, including the randomly collected data in \Cref{sec:collecting_statistics} and the solver selection experiment \Cref{sec:solver_selection_experiment}, and rerun the experiments with it. This mimics the linear solver expert with a strong prior knowledge of the considered models. We therefore denote this case as ``expert''. The incremental updates are still enabled for the expert case, as is the capability to switch the solver configuration during the simulation, but the initial exploration is disabled. As before, we repeat this experiment five times with the simulations in the shuffled order. 

\Cref{fig:runtime_expert} compares the run times of the expert case and the solver selection case, described in \Cref{sec:solver_selection_experiment}. The plots show cumulative run times, including the time to construct and solve linear systems and the machine learning overhead. Five lines of the same color correspond to five repeats of each sequence, and bold lines correspond to the mean across five runs. The time spent on failed attempts to solve linear systems is also included, which is reflected in jumps at the beginning of the sequences for the solver selection case. For the expert case, all solve attempts were successful. The solver selection case run time includes 64 linear systems corresponding to the initial exploration. 

Bad decisions during and after the initial exploration have a stronger impact in Sequence A for two reasons. First, the solver space includes 35\% configurations that always fail, versus 8\% for Sequence B, according to \Cref{fig:statistics_random_runs}. This indicates that 64 configuration in the initial exploration is not enough for Sequence A to avoid bad decisions, and it needs to make more ``mistakes'', which cause large jumps on the left part of the figure. However, after $\sim$200 linear systems, the solver selection algorithm stabilizes and starts to avoid bad decisions, and the lines become parallel to the ``expert'' case.
The other reason is that each run of Sequence B comes with almost twice as many linear systems, making the cumulative overhead of bad decisions in the left side of the figure negligible.

For Sequence B, the expert case performs on average marginally better than the solver selection case. This illustrates that it is difficult to improve on
initial performance data collection, and inexpensive tuning on the fly results in a nearly optimal performance.
We observe that in Sequence B, some runs of the solver selection case overtake the expert case. This is attributed to the increased machine-learning overhead of the expert case, which needs to retrain the machine learning pipeline based on all its data after every 64 linear systems.

\section{Conclusion}
\label{sec:conclusion}
This paper addresses the problem of preconditioned linear solver selection and tuning for multiphysics simulations, considering two time-dependent nonlinear model problems: (A) coupled flow and heat transfer in porous media; and (B) thermo-poromechanics in porous media with fractures, governed by frictional contact mechanics. The problem poses three main challenges: (i) the large number of preconditioner sub-algorithm combinations results in thousands of possible solver configurations; (ii) solver performance depends on the solution at each time step, as nonlinear couplings influence matrix coefficients during linearization; and (iii) the solver selection algorithm must operate in a setting where performance data from previous simulations quickly becomes irrelevant, so it is not possible to collect a large performance dataset.

We propose the preconditioned linear solver selection algorithm based on a pipeline of two machine learning models: a classifier that screens the solver configurations and discards likely failures, and a regression model that predicts the linear solver run time. Based on the pipeline, we select a solver configuration for each new linear system. The machine learning pipeline is updated with feedback during the simulation, resulting in a continuously refined selection policy.

The numerical experiments include two sequences of simulations of the model problems, where each sequence consists of several experiments of the same model with slightly different simulation setups. We apply the solver selection algorithm to them and show that (i) it is capable of avoiding failures, and (ii) it manages to select the options within the best 3\% of all considered configurations.
These results indicate that the proposed solver selection algorithm can be used in practice, resulting in robust and efficient simulations. The cost of this result is solving a few linear systems with random linear solver configurations in the beginning of a simulation, and a negligible machine learning overhead, relatively to the linear solver run time.

Finally, we compare the performance of the simulations using the solver selection algorithm with that of simulations guided by a machine learning model that has access to all previously collected data. The latter serves as an ``expert'' with strong prior knowledge of linear solver performance. Results show that after the necessary exploration, the ``expert'' case performs only marginally better than the proposed solver selection approach. This suggests that the solver selection algorithm effectively addresses the challenge of solver selection and tuning, which is particularly valuable for simulation engineers and researchers, given that human expertise in tuning linear solvers is not always readily available.

\section*{Acknowledgment}
This project has received funding from the VISTA program, The Norwegian Academy of Science and Letters and Equinor and from the European Research Council (ERC) under the European Union’s Horizon 2020 research and innovation program (grant agreement No 101002507).

\bibliographystyle{unsrt}
\bibliography{references}

\end{document}